\documentclass[a4paper,11pt]{article}

\usepackage{natbib}
\usepackage[english]{babel}
\usepackage{amsmath}
\usepackage{amssymb}
\usepackage{fancyhdr}
\usepackage[latin1]{inputenc}
\usepackage{stmaryrd}
\usepackage{amsthm}
\usepackage{geometry}
\usepackage{enumerate}
\usepackage{hyperref}
\usepackage{graphicx}

\newtheorem{thm}{Theorem}[subsection]
\newtheorem{prop}{Proposition}[subsection]
\newtheorem{lm}{Lemma}[subsection]
\newcommand{\e}{\varepsilon}
\newcommand{\E}{\mathbb{E}}
\newcommand{\dd}{\mathrm{d}}
\newcommand{\OO}{\mathrm{O}}
\newcommand{\R}{\mathrm{Re}}

\newcommand{\tr}{\mathrm{trace}}

\makeatletter
\renewcommand{\thefigure}{\ifnum \c@section>\z@ \thesection.\fi
 \@arabic\c@figure}
\@addtoreset{figure}{section}
\makeatother

\newcommand{\auteur}[4]{ \author{ { #1}\footnote{#2} \vspace{0.5cm} 
 \\  
 \hspace{-2.0cm} \small \noindent\begin{tabular}{l}   #3 \\ #4\end{tabular}}}

\auteur{Fanny Godet}{fanny.godet@math.univ-nantes.fr}{Laboratoire de Mathématiques Jean Leray, UMR CNRS 6629  } { Université de Nantes
2 rue de la Houssinière - BP 92208 
F-44322 Nantes Cedex 3 }
\title{Linear Prediction of Long-Range Dependent Time Series}

\begin{document}
\maketitle
\begin{abstract}
We present two approaches for next step linear prediction of long memory time series. The first is based on the truncation of the Wiener-Kolmogorov predictor by restricting the observations to the last $k$ terms, which are the only available values in practice. Part of the mean squared prediction error comes from the truncation, and another part comes from the parametric estimation of the parameters of the predictor. By contrast, the second approach is non-parametric. An AR($k$) model is fitted to the long memory time series and we study the error made with this misspecified model.
\end{abstract}


\par ARMA (autoregressive moving-average) processes are often called short-memory processes because their covariances decay rapidly (i.e. their covariance decay exponentially). By contrast, a long-memory process is characterised by the following feature: the autocovariance function $\sigma$  decays more slowly i.e. it is not absolutely summable. They are so-named because of the strong association between observations widely separated in time. The long-memory time series models have attracted much attention lately and there is now a growing realisation that time series possessing long-memory characteristics arise in subject areas as diverse as Economics, Geophysics, Hydrology or telecom traffic (see, e.g., \cite{mandelbrot} and \cite{granger}). Although there exists substantial literature on the prediction of short-memory processes(see \cite{bhansali} for the univariate case or \cite{lewis} for the multivariate case), there are less results for long-memory time series. In this paper, we consider the question of the prediction of the latter.
\par More precisely, we will compare two prediction methods for long-memory process. Our goal is a linear predictor $\widetilde{X}_{k+1}$ from observed values which is optimal in the sense that it minimizes the mean-squared error $\mathbb{E}\left[  \left( X_{k+1}-\widetilde{X}_{k+1}\right) ^2\right] $. This paper is organized as follows. First we will introduce our model and our main assumptions. Then in section \ref{meilleurpredicteurlineairetronque}, we study the best linear predictor i.e. the Wiener-Kolmogorov predictor proposed by \cite{whittle} and by \cite{tronc} for the long-memory time series. In practice, only the last $k$ values of the process are available. Therefore we need to truncate the infinite series which defines the predictor and derive the asymptotic behaviour as $k \rightarrow + \infty$ of the mean-squared error. Then we propose an estimator of the coefficients of the infinite autoregressive representation based on a realisation of length $T$. Under the simplifying assumption that the series used for estimation and the series used for prediction are generated from two independent process which have the same stochastic structure, we obtain an approximation of the mean-squared prediction error when $T \rightarrow + \infty$ and then $k \rightarrow + \infty$. 
\par In Section \ref{sectionarp}, we discuss the asymptotic properties of the forecast error if we fit a misspecified AR($k$) model to a long-memory time series. This approach has been proposed by \cite{arp} for fractional noise series F($d$). His simulations show that high-order AR m-models forecast fractional integrated noise very well. In that case we also study the consequences of the estimation of the forecast coefficients. Therefore we shall rewrite the heuristic proof of Theorem 1 of \cite{arp} and develop a generalization of this result to a larger class of long-memory models. We conclude by comparing our asymptotic approximation for the global prediction error of long-memory processes and that of \cite{berk} and \cite{bhansali} in the case of short memory time series. Subsidiary proofs are given in the Appendix.

\section{Model}
\label{intro}
\par Let $(X_n)_{n \in \mathbb{Z}}$ be a discrete-time (weakly) stationary process in $\mathrm{L}^2$ with mean 0 and $\sigma$ its autocovariance function. We assume that the process $(X_n)_{n \in \mathbb{Z}}$ is a long-memory process i.e.:
\begin{displaymath}
\sum_{k=-\infty}^{\infty} |\sigma(k)|=\infty.
\end{displaymath}
The process $(X_n)_{n \in \mathbb{Z}}$ admits an infinite moving average representation as follows:
\begin{equation}
 X_n=\sum_{j=0}^{\infty} b_j \e_{n-j} \label{defbj}
\end{equation}
where $(\varepsilon_n)_{n\in \mathbb{Z}}$ is a white-noise series consisting of uncorrelated random variables, each with mean 0 and variance $\sigma^2_{\varepsilon}$ and $(b_j)_{j \in \mathbb{N}}$ are square-summable.
 We shall further assume that $(X_n)_{n \in \mathbb{Z}}$ admits an infinite autoregressive representation:
\begin{equation}
\e_n=\sum_{j=0}^{\infty} a_j X_{n-j}, \label{defaj}
\end{equation}
where the $(a_j)_{j \in \mathbb{N}}$ are absolutely summable.
We assume also that $(a_j)_{j \in \mathbb{N}}$ and  $(b_j)_{j \in \mathbb{N}}$, occurring respectively in \eqref{defaj} and \eqref{defbj}, satisfy the following conditions for all $\delta>0$:
\begin{eqnarray}
|a_j|&\leq&C_1  j^{-d-1+\delta} \label{aj} \\
|b_j|&\leq& C_2 j^{d-1+\delta} .\label{bj}
\end{eqnarray}
where $C_1$ and $C_2$ are constants and $d$ is a parameter verifying $d \in ]0,1/2[$. For example, a FARIMA process $(X_n)_{n \in \mathbb{Z}}$ is the stationary solution to the difference equations:
\begin{displaymath}
\phi(B) (1-B)^d X_n = \theta(B) \e_n
\end{displaymath}
where $(\e_n)_{n \in \mathbb{Z}}$ is a white noise series, $B$ is the backward shift operator and $\phi$ et $\theta$ are polynomials with no zeroes on the unit disk. Its coefficients verify equations \eqref{aj} and \eqref{bj}. In particular, if $\phi=\theta=1$ then the process $(X_n)_{n \in \mathbb{Z}}$ is called \textit{fractionally integrated noise} and denoted F($d$). More generally, series like:
\begin{eqnarray}
|a_j|&\underset{+\infty}\sim& L(j^{-1}) j^{-d-1}  \nonumber \\
|b_j|&\underset{+\infty}\sim& L'(j^{-1}) j^{d-1} \nonumber 
\end{eqnarray}
where $L$ and $L'$ are slowly varying functions and therefore verify conditions \eqref{aj} and \eqref {bj}. 
A positive $L$ will be called a slowly varying function in the sense of \cite{zygmund} if, for any $\delta>0$, $x \mapsto x^{-\delta}L(x)$ is decreasing and $x \mapsto x^{\delta}L(x)$ is increasing. 
\par  The condition \eqref{bj} implies that the autocovariance function $\sigma$ of the process $(X_n)_{n \in \mathbb{Z}}$ verifies:
\begin{equation}
\forall \delta>0, \exists C_3 \in \mathbb{R} ,\quad \vert \sigma(j)\vert \leq  C_3j^{2d-1+\delta} .\label{sigma}
\end{equation}
Since, if $\delta<\frac{1-2d}{2}$:
\begin{eqnarray*}
\sigma(k)&=& \sum_{j=0}^{+\infty} b_j b_{j+k} \\
\vert \sigma(k) \vert &\leq&  \sum_{j=0}^{+\infty}\vert b_j b_{j+k}\vert\\
&\leq& C_2^2\sum_{j=0}^{+\infty} j^{d-1+\delta} (k+j)^{d-1+\delta} \\
&\leq& C_2^2 \int_{-1}^{+\infty} j^{d-1+\delta} (k+j)^{d-1+\delta} \dd j \\
&\leq& C_2^2k^{2d-1+2\delta}\int_{-1}^{+\infty}j^{d-1+\delta} (1+j)^{d-1+\delta} \dd j\\
&\leq& C_3k^{2d-1+2\delta}
\end{eqnarray*}
Notice that it suffices to prove \eqref{sigma} for $\delta$ near 0 in order to verify  \eqref{sigma} for $\delta >0$ arbitrarily chosen. More accurately, \cite{inoueregularly} has proved than if:
\begin{displaymath}
b_j \sim L\left( j^{-1}\right) j^{d-1} 
\end{displaymath}
then 
\begin{displaymath}
\sigma(j) \sim j^{2d-1} \left[ L\left( j^{-1}\right)\right]^2 \beta(1-2d,d)
\end{displaymath}
where $L$ is a slowly varying function and $\beta$ is the beta function. The converse is not true, we must have more assumptions about the series $(b_j)_{j\in \mathbb{N}}$ in order to get an asymptotic equivalent for $(\sigma(j))_{j\in \mathbb{N}}$ (see \cite{inouepartialautocorr}).

\section{Wiener-Kolmogorov Prediction Theory}
\label{meilleurpredicteurlineairetronque}
\par The aim of this part is to compute the best linear one-step predictor (with minimum mean-square distance from the true random variable) knowing all the past $\{X_{k+1-j},j \leqslant 1\}$. Our predictor is therefore an infinite linear combination of the infinite past:
\begin{equation}
\widetilde{X_k}(1)=\sum_{j=0}^{\infty} \lambda(j) X_{k-j} \nonumber 
\end{equation}
 where $(\lambda(j))_{j \in \mathbb{N}}$ are chosen to ensure that the mean squared prediction error:
\begin{displaymath}
\mathbb{E}\big[ \big(\widetilde{X_k}(1)-X_{k+1}\big)²\big]
\end{displaymath}
is as small as possible. Following \cite{whittle}, and in view of the moving average representation of $(X_n)_{n \in \mathbb{Z}}$, we may rewrite our predictor $\widetilde{X_k}(1)$ as:
\begin{displaymath}
\widetilde{X_k}(1)=\sum_{j=0}^{\infty}\phi(j)\varepsilon_{k-j}.
\end{displaymath}
where $(\phi(j))_{j \in \mathbb{N}}$ depends only on $(\lambda(j))_{j \in \mathbb{N}}$ and $(a_j)_{j \in \mathbb{N}}$ defined in \eqref{defaj}.
From the infinite moving average representation of $(X_n)_{n \in \mathbb{Z}}$ given below in \eqref{defbj}, we can rewrite the mean-squared prediction error as:
\begin{eqnarray*}
\mathbb{E}\big[ \big(\widetilde{X_k}(1)-X_{k+1}\big)²\big]
&=&\mathbb{E}\left[ \left( \sum_{j=0}^{\infty}\phi(j)\varepsilon_{k-j}-\sum_{j=0}^{\infty} b(j)\varepsilon_{k+1-j}\right) ^2\right]  \\  
&=&\mathbb{E}\left[ \left(\varepsilon_{k+1}-\sum_{j=0}^{\infty}\left( \phi(j)-b(j+1)\right) \varepsilon_{k-j}\right) ^2\right] \\
&=&\left( 1+\sum_{j=0}^{\infty}\big( b_{j+1}-\phi(j)\big)²\right) \sigma_{\e}²
\end{eqnarray*}
since the random variables $(\e_n)_{n \in \mathbb{Z}}$ are uncorrelated with variance $\sigma_{\e}²$. The smallest mean-squared prediction error is obtained when setting $ \phi(j)=b_{j+1}$ for $j\geq 0$. 
 \par The smallest prediction error of $(X_n)_{n \in \mathbb{Z}}$ is $\sigma_{\e}²$ within the class of linear predictors.
Furthermore, if 
\begin{displaymath}
A(z)=\sum_{j=0}^{+\infty} a_j z^j,
\end{displaymath}
denotes the characteristic polynomial of the  $(a(j))_{j \in \mathbb{Z}}$ and
\begin{displaymath}
B(z)=\sum_{j=0}^{+\infty} b_j z^j,
\end{displaymath}
that of the $(a(j))_{j \in \mathbb{Z}}$, then in view of the identity, $A(z)=B(z)^{-1},\,\vert z \vert \leq 1$, we may write:
\begin{equation}
\widetilde{X_k}(1)=-\sum_{j=1}^{\infty}a_jX_{k+1-j} .\label{serieinf}
\end{equation}

\subsection{Mean Squared Prediction Error when the Predictor is Truncated}
\label{erreurtronc}
\par In practice, we only know a finite part of the past, the one which we have observed. So the predictor should only depend on the observations. Assume that we only know the set $\{X_0, \ldots, X_k\}$ and that we replace the unknown values by 0, then we have the following new predictor: 
\begin{equation}
\widetilde{X'_k}(1)=-\sum_{j=1}^{k}a_jX_{k+1-j} .\label{pred1}
\end{equation}
It is equivalent to say that we have truncated the infinite series \eqref{serieinf} to $k$ terms.
The following proposition provides us the asymptotic properties of the mean squared prediction error as a function of $k$.

\begin{prop} Let $(X_n)_{n \in \mathbb{Z}}$ be a linear stationary process defined by \eqref{defbj}, \eqref{defaj} and possessing the features \eqref{aj} and \eqref{bj}. We can approximate the mean-squared prediction error of $\widetilde{X'_k}(1)$ by:
\begin{displaymath}
\forall \delta>0, \quad \mathbb{E}\big( \big[X_{k+1}-\widetilde{X'_k}(1)\big]² \big)  =\sigma_{\e}²+\OO(k^{-1+\delta}).
\end{displaymath}
Furthermore, this rate of convergence $\OO(k^{-1})$ is optimal since for fractionally integrated noise, we have the following asymptotic equivalent:
\begin{displaymath}
\mathbb{E}\big( \big[X_{k+1}-\widetilde{X'_k}(1)\big]² \big)  =\sigma_{\e}²+C k^{-1}+\mathrm{o}\left(k^{-1} \right) .
\end{displaymath}

\end{prop}
We note that the prediction error is the sum of $\sigma_{\e}²$, the error of Wiener-Kolmogorov model and the error due to the truncation to $k$ terms which is bounded by $\OO(k^{-1+\delta})$ for all $\delta>0$.
\begin{proof}
\begin{eqnarray}
X_{k+1}-\widetilde{X'_k}(1) &=& X_{k+1}-\widetilde{X_k}(1)+\widetilde{X_k}(1)-\widetilde{X'_k}(1) \nonumber \\
&=& X_{k+1} -\sum_{j=0}^{+\infty}b_{j+1} \e_{k-j} -\sum_{j=k+1}^{+\infty}a_j X_{k+1-j} \nonumber \\
&=& \e_{k+1}-\sum_{j=k+1}^{+\infty}a_j X_{k+1-j}  \label{orthogonalite}.
\end{eqnarray}
The two parts of the sum \eqref{orthogonalite} are orthogonal for the inner product associated with the mean square norm. Consequently:
\begin{displaymath}
\mathbb{E}\big( \big[X_{k+1}-\widetilde{X'_k}(1)\big]² \big) 
 = \sigma_{\e}²+ \sum_{j=k+1}^{\infty} \sum_{l=k+1}^{\infty} a_j a_l \sigma(l-j)  .
\end{displaymath}
For the second term of the sum we have:
\begin{eqnarray}
\bigg \vert \sum_{j=k+1}^{+\infty} \sum_{l=k+1}^{+\infty} a_j a_l \sigma(l-j) \bigg \vert&=&\bigg \vert 2 \sum_{j=k+1}^{+\infty} a_j \sum_{l=j+1}^{+\infty}a_l \sigma(l-j) + \sum_{j=k+1}^{+\infty}a_j^2 \sigma(0)\bigg \vert \nonumber \\
&\leq & 2 \sum_{j=k+1}^{+\infty}|a_j|\left|a_{j+1}\right|| \sigma(1)|+ \sum_{j=k+1}^{+\infty}a_j^2 \sigma(0) \nonumber\\
&& +2\sum_{j=k+1}^{+\infty} |a_j|\sum_{l=j+2}^{+\infty}|a_l|| \sigma(l-j)|\nonumber
\end{eqnarray}
from the triangle inequality, it follows that:
\begin{eqnarray}
&&\bigg \vert \sum_{j=k+1}^{+\infty} \sum_{l=k+1}^{+\infty} a_j a_l \sigma(l-j) \bigg \vert \nonumber \\
&\leq &
C_1^2C_3\left( 2\sum_{j=k+1}^{+\infty}j^{-d-1+\delta}(j+1)^{-d-1+\delta}+\sum_{j=k+1}^{+\infty}\left( j^{-d-1+\delta}\right)^2 \right) \label{sommesimple}\\
&+&2C_1^2C_3\sum_{j=k+1}^{+\infty}j^{-d-1+\delta} \sum_{l=j+2}^{+\infty}l^{-d-1+\delta}\vert l-j\vert^{2d-1+\delta} \label{sommedouble} 
\end{eqnarray}
 for all $\delta>0$ from inequalities \eqref{aj} and \eqref{sigma}. Assume now that $\delta <1/2-d$. For the terms \eqref{sommesimple}, since $j \mapsto j^{-d-1+\delta} (j+1)^{-d-1+\delta}$ is a positive and decreasing function on $\mathbb{R}^+$, we have the following approximations:
 \begin{eqnarray*}
 2C_1^2C_3\sum_{j=k+1}^{+\infty}j^{-d-1+\delta}(j+1)^{-d-1+\delta} & \sim &2C_1^2C_3 \int_k^{+\infty}j^{-d-1+\delta}(j+1)^{-d-1+\delta} \dd j \\
  & \sim &\frac{2C_1^2C_3}{1+2d-2\delta} k^{-2d-1+2\delta} 
 \end{eqnarray*}
Since the function $j \mapsto \left( j^{-d-1+\delta}\right) ^2$ is also positive and decreasing, we can establish in a similar way that:
 \begin{eqnarray*}
 C_1^2C_3\sum_{j=k+1}^{+\infty}\left( j^{-d-1+\delta}\right) ^2 & \sim&C_1^2C_3\int_k^{+\infty}\left( j^{-d-1+\delta}\right) ^2\dd j \\
 & \sim &\frac{C_1^2C_3}{1+2d-2\delta} k^{-2d-1+2\delta}. 
 \end{eqnarray*}

For the infinite double series \eqref{sommedouble}, we will similarly compare the series with an integral. In the next Lemma, we establish the necessary result for this comparison:
\begin{lm} \label{sommedoublecomp}
Let $g$ the function $(l,j) \mapsto j^{-d-1+\delta} \, l^{-d-1+\delta}\,\vert l-j\vert^{2d-1+\delta}$. Let $m$ and $n$ be two positive integers. We assume that $\delta<1-2d $ and $m \geq \frac{\delta -d-1}{\delta+2d-1}$ for all $\delta \in \left] 0,\frac{\delta -d-1}{\delta+2d-1}\right[ $. 
We will call $A_{n,m}$ the square $[n,n+1]\times[m,m+1]$. If $n \geq m+1$ then
\begin{displaymath}
\int_{A_{n,m}}g(l,j)\, \dd j \,\dd l \geq g(n+1,m).
\end{displaymath}
\end{lm}
\begin{proof}
see the appendix \ref{annexe1}
\end{proof}

Assume now that $\delta<1-2d $ without loss of generality.
Thanks to the previous Lemma and the asymptotic equivalents of \eqref{sommesimple}, there exists $K \in \mathbb{N}$ such that if $k>K$:
\begin{eqnarray}
\bigg \vert \sum_{j=k+1}^{+\infty} \sum_{l=k+1}^{+\infty} a_j a_l \sigma(l-j) \bigg \vert&\leq & C \int_{k+1}^{+\infty} j^{-d-1+\delta}\left[  \int_{j}^{+\infty} l^{-d-1+\delta} (l-j)^{2d-1+\delta} \dd l \right]  \dd j +\OO\left( k ^{-2d-1+2\delta}\right) \nonumber
\end{eqnarray}
In the integral over $l$ by using the substitution $jl'=l$, we obtain:
\begin{displaymath}
\bigg \vert \sum_{j=k+1}^{+\infty} \sum_{l=k+1}^{+\infty} a_j a_l \sigma(l-j) \bigg \vert\leq  C'\int_{k+1}^{+\infty}j^{-2+3\delta} \int_{1}^{+\infty} l^{-d-1+\delta} (l-1)^{2d-1+\delta} \dd l \dd j +\OO\left( k ^{-2d-1}\right).  
\end{displaymath} Since if $\delta<(1-d)/2$
\begin{displaymath}
\int_1^{+\infty} l^{-d-1+\delta}(l-1)^{2d-1+\delta} \dd l < +\infty,
\end{displaymath}
it follows:
\begin{eqnarray}
\bigg \vert \sum_{j=k+1}^{+\infty} \sum_{l=k+1}^{+\infty} a_j a_l \sigma(l-j) \bigg \vert&\leq & \OO\left(k^{-1+3\delta} \right)+\OO\left( k ^{-2d-1}\right) \nonumber\\
&\leq & \OO\left(k^{-1+3\delta} \right).
\end{eqnarray}

If $\delta>0$, $\delta<1-2d $ and $\delta<(1-d)/2$, we have:
\begin{displaymath}
\bigg \vert \sum_{j=k+1}^{+\infty} \sum_{l=k+1}^{+\infty} a_j a_l \sigma(l-j) \bigg \vert = \OO\left(k^{-1+3\delta} \right).
\end{displaymath}
Notice that if the equality is true under the assumptions $\delta>0$, $\delta<1-2d $ and $\delta<(1-d)/2$, it is also true for any $\delta>0$. Therefore we have proven the first part of the theorem. \\
 We prove now that there exists long-memory processes whose prediction error attains the rate of convergence $k^{-1}$. Assume now that $(X_n)_{n \in \mathbb{Z}}$ is fractionally integrated noise F($d$), which is the stationary solution of the difference equation:
\begin{equation}
X_n=(1-B)^{-d} \e_n \label{Fd}
\end{equation}
with $B$ the usual backward shift operator, $(\e_n)_{n \in \mathbb{Z}}$ is a white-noise series and $d\in \left]0, 1/2\right[ $ (see for example \cite{bd}). We can compute the coefficients and obtain that:
\begin{displaymath}
\forall j>0,\quad a_j=\frac{\Gamma(j-d)}{\Gamma(j+1)\Gamma(-d)} \textrm{ and } \forall j\geq 0,\quad \sigma(j)=\frac{(-1)^j\Gamma(1-2d)}{\Gamma(j-d+1)\Gamma(1-j-d)} \sigma_\e^2
\end{displaymath}
then we have:
\begin{displaymath}
\forall j>0,\quad a_j<0\textrm{ and } \forall j\geq 0,\quad \sigma(j)>0
\end{displaymath}
and 
\begin{displaymath}
a_j \sim \frac{j^{-d-1}}{\Gamma(-d)}\textrm{ and }\sigma(j)\sim \frac{j^{2d-1}\Gamma(1-2d)}{\Gamma(d)\Gamma(1-d)} \quad \textrm{when } j \rightarrow \infty. 
\end{displaymath}

In this particular case, we can estimate the prediction error more precisely:
\begin{eqnarray}
\sum_{k+1}^{+\infty} \sum_{k+1}^{+\infty} a_j a_l \sigma(l-j)&=&\sum_{k+1}^{+\infty} |a_j|\sum_{j+1}^{+\infty}|a_l|| \sigma(l-j)| + \sum_{k+1}^{+\infty}a_j^2 \sigma(0) \nonumber\\
&\sim&\frac{\Gamma(1-2d)}{\Gamma(-d)^2\Gamma(d)\Gamma(1-d)}\int_{k+1}^{+\infty}j^{-2} \int_{1/j+1}^{+\infty} l^{-d-1} (l-1)^{2d-1} \dd l \dd j +\OO\left( k ^{-2d-1}\right) \nonumber\\
\sum_{k+1}^{+\infty} \sum_{k+1}^{+\infty} a_j a_l \sigma(l-j) &\sim& \frac{\Gamma(1-2d)\Gamma(2d)}{\Gamma(-d)^2\Gamma(d)\Gamma(1+d)} k^{-1} \label{premest}
\end{eqnarray}
The asymptotic bound $\OO(k^{-1})$ is therefore as small as possible.
\end{proof}

\par In the specific case of fractionally integrated noise, we may write the prediction error as:
\begin{displaymath}
\mathbb{E}\big( \big[X_{k+1}-\widetilde{X'_k}(1)\big]² \big) 
 = \sigma_{\e}²+C(d) k^{-1}+\mathrm{o}\left( k^{-1}\right) 
\end{displaymath}
and we can express $C(d)$ as a function of $d$:
\begin{equation}
C(d)=\frac{\Gamma(1-2d)\Gamma(2d)}{\Gamma(-d)^2\Gamma(d)\Gamma(1+d)}\label{C}.
\end{equation}
 It is easy to prove that $C(d) \rightarrow + \infty$ as $d \rightarrow 1/2$ and we may write the following asymptotic equivalent as $d \rightarrow 1/2$:
 \begin{equation}
 C(d) \sim \frac{1}{(1-2d)\Gamma(-1/2)^2\Gamma(1/2)\Gamma(3/2)}.\label{eqcd}
 \end{equation}
As $d \rightarrow 0 $, $C(d) \rightarrow 0$ and we have the following equivalent as $d \rightarrow 0 $:
 \begin{displaymath}
 C(d) \sim d^2.
 \end{displaymath}

\begin{figure}[h]
\begin{center}
\caption{The Constant $C(d)$, $d \in [0,1/2[$, defined in \eqref{C}}
\label{figureC}
\includegraphics[width=10cm]{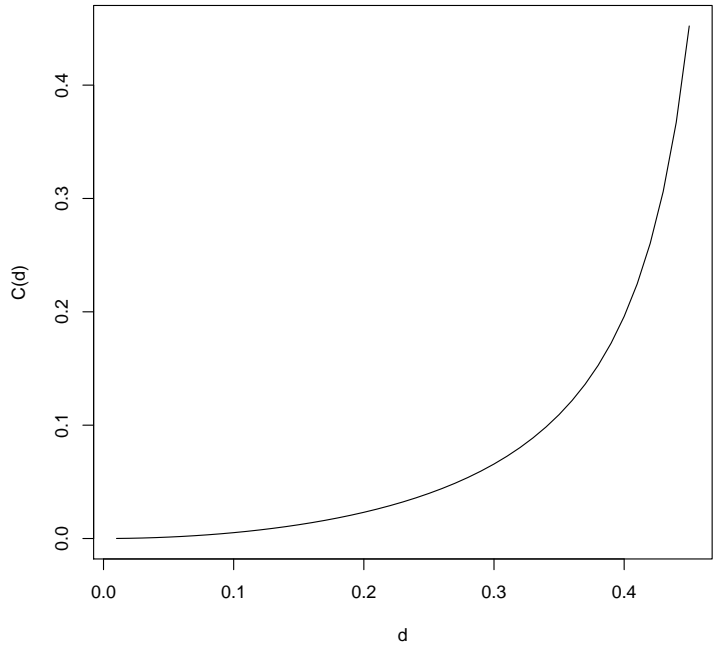}
\end{center}
\end{figure}
\par As the figure \ref{figureC} suggests and the asymptotic equivalent given in \eqref{eqcd} proves, the mean-squared error tends to $+\infty$ as $d \rightarrow 1/2$. By contrast, the constant $C(d)$ takes small values for d in a large interval of $[0,1/2[$.Although the rate of convergence has a constant order $k^{-1}$, the forecast error is bigger when $d \rightarrow 1/2$. This result is not surprising since the correlation between the random variable, which we want to predict, and the random variables, which we take equal to 0, increases when $d \rightarrow 1/2$.

\subsection{Estimates of Forecast Coefficients and the Associated Mean Square Error}
\label{parametrique}

We will now estimate the mean-squared error between the predictor $\widetilde{X'_k}(1)$ defined on \eqref{pred1} and the predictor $\widetilde{X'_{T,k}}(1)$ defined as:
\begin{displaymath}
\widetilde{X'_{T,k}}(1):=-\sum_{j=1}^{k}\widehat{a}_jX_{k+1-j} 
\end{displaymath}
where $\widehat{a}_j$ are estimates of $a_j$ computed using a length $T$ realisation of the process.
More precisely, we consider a parametric approach and we assume that:
\begin{displaymath}
a_j=a_j(\theta) \textrm{ with } \theta \textrm{ an unknown vector in } \Theta
\end{displaymath}
where $\Theta$ is a compact subset of $\mathbb{R}^p$.
 Assume that the process $(Y_n)_{n \in \mathbb{Z}}$ is Gaussian. Let $\theta_0$ be the true value of the parameter. We assume the realisation $(Y_n)_{1 \leq n \leq T}$ to be known. We estimate the $(a_j)_{1 \leq j \leq k}$ by $\widehat{a}_j:=a_j(\widehat{\theta}_T)$ where $\widehat{\theta}_T$ is an estimate of $\theta_0$, for example the Whittle estimate. In order to use the Whittle estimate and follow the approach suggested in \cite{foxtaqqu}, we assume from now on that all the processes in the parametric class have a spectral density denoted by $f(.,\theta)$.

We define the Whittle estimate by (see \cite{foxtaqqu}):
\begin{equation}
\hat{\theta}_T=\underset{\theta \in \Theta}{\mathrm{argmin}}\left[ \frac{1}{2 \pi} \int_{-\pi}^{\pi}\left[ f(\lambda,\theta)\right]^{-1} I_T(\lambda) \dd \lambda \right] \label{Whittle}
\end{equation}
where $I_T$ is the periodogram:
\begin{displaymath}
I_T(\lambda)=\frac{\vert\sum_{j=1}^T \mathrm{e}^{ij \lambda}(Y_j-\overline{Y_T})\vert^2}{2 \pi T}.
\end{displaymath}
\par Before we state the theorem, we will give assumptions on the regularity of the spectral densities in our parametric class. Under those standard conditions, the estimated vector converges to the true parameter if the process is a Gaussian long-memory time series (see  \cite{foxtaqqu}).\\We will refer to the following assumptions. 
\\ We say that $f(x,\theta)$ satisfies conditions A0-A6 if there exists $0<\alpha(\theta)<1$ such that for each $\delta >0$,

\vspace{0.3 cm}
\noindent
A0. $f(\lambda,\theta_0)=|\lambda|^{2\alpha(\theta_0)}L(\lambda,\theta_0)$ with $L(.,\theta_0)$ bounded. $L(.,\theta_0)$ is differentiable at 0 and $L(.,\theta_0)\neq 0$. 

\vspace{0.2 cm}

\noindent
A1. $\theta \mapsto \int_{-\pi}^{\pi}f(\theta,\lambda) \dd \lambda < + \infty$ can be twice differentiated under the integral sign.
\vspace{0.2 cm}

\noindent
A2. $f(\theta,\lambda)$ is continuous at all $(\theta,\lambda)$, $\lambda\neq0$, $f^{-1}(\theta,\lambda)$ is continuous at all $(\theta,\lambda)$ and, 
\begin{displaymath}
f(\theta,\lambda)=\OO(\vert \lambda\vert^{-\alpha(\theta)-\delta }) \quad \textrm{ as } \lambda \rightarrow 0.
\end{displaymath}

\noindent
A3. $(\partial/\partial \theta_j) f^{-1}(\theta,\lambda)$ and $(\partial^2/\partial \theta_j\partial \theta_l) f^{-1}(\theta,\lambda)$ are continuous at all $(\theta,\lambda)$, 
\begin{displaymath}
\forall 1 \leq j \leq p, \quad \frac{\partial}{\partial \theta_j}f^{-1}(\theta,\lambda)=\OO(\vert \lambda\vert^{\alpha(\theta)-\delta }) \quad \textrm{ as } \lambda \rightarrow 0
\end{displaymath}
and
\begin{displaymath}
\forall 1 \leq j,l \leq p\quad \frac{\partial^2}{\partial \theta_j\partial \theta_l}f^{-1}(\theta,\lambda)=\OO(\vert \lambda\vert^{\alpha(\theta)-\delta }) \quad \textrm{ as } \lambda \rightarrow 0.
\end{displaymath}

\vspace*{-0.5em}

\noindent
A4. $(\partial/\partial \lambda) f(\theta,\lambda)$ is continuous at all $(\theta,\lambda)$, $\lambda\neq0$, and
\begin{displaymath}
\frac{\partial}{\partial \lambda}f(\theta,\lambda)=\OO(\vert \lambda\vert^{-\alpha(\theta)-1-\delta }) \quad \textrm{ as } \lambda \rightarrow 0
\end{displaymath}

\noindent
A5. $(\partial^2/\partial \theta_j\partial \lambda) f^{-1}(\theta,\lambda)$ are continuous at all $(\theta,\lambda)$, $\lambda\neq0$, and 
\begin{displaymath}
\forall 1 \leq j \leq p, \quad \frac{\partial^2}{\partial \theta_j\partial \lambda}f^{-1}(\theta,\lambda)=\OO(\vert \lambda\vert^{\alpha(\theta)-1-\delta }) \quad \textrm{ as } \lambda \rightarrow 0.
\end{displaymath}

\noindent
A6. $(\partial^3/\partial \theta_j\partial^2 \lambda) f^{-1}(\theta,\lambda)$ are continuous at all $(\theta,\lambda)$, $\lambda\neq0$, and
\begin{displaymath}
\forall 1 \leq j \leq p, \quad \frac{\partial^3}{\partial \theta_j\partial^2 \lambda}f^{-1}(\theta,\lambda)=\OO(\vert \lambda\vert^{\alpha(\theta)-2-\delta })
\end{displaymath}

We can now express the asymptotic behavior of the mean-squared prediction error due to the estimation of the forecast coefficients. We assume in this Section that the process is Gaussian. 
Let $(X_j)_{j \in \mathbb{Z}}$ be a stochastic process, which verifies the assumptions of section \ref{intro}, and let $(Y_j)_{j \in \mathbb{Z}}$ be a process which is independent of $(X_j)_{j \in \mathbb{Z}}$, but has the same stochastic structure. We want to predict $X_{k+1}$ knowing $(X_j)_{j \in \llbracket 1 ,k \rrbracket}$ and we assume that the parameter $\theta$ and so the forecast coefficients are estimated based on a realisation  $(Y_j)_{j \in \llbracket 1 ,T \rrbracket}$.

\begin{thm}
Let $(X_n)_{n \in \mathbb{Z}}$ be a stationary Gaussian long-memory sequence with mean 0 and spectral density $f(\theta,\lambda)$ and $\theta \in \Theta$ is an unknown parameter. The set $\Theta$ is assumed to be compact. We assume also that $\theta_0$ is in the interior of $\Theta$ and that $\forall \theta \in \mathring{\Theta}$, the conditions A1-A6 hold.
Moreover we assume that each process $(Z_n)_{n \in \mathbb{Z}}$ in our parametric class with $\theta \in \Theta$ admits an autoregressive representation:
\begin{displaymath}
\e_n=\sum_{j=0}^{\infty} a_j(\theta) Z_{n-j}
\end{displaymath}
where $(\e_n)_{n \in \mathbb{Z}}$ is a Gaussian white noise.
Let $\theta_0$ be the true value of the parameter and assume that $\theta_0 \in \mathring{\Theta}$. Assume also that $f(\theta_0,\lambda)$ verifies A0 and that for any $j \in \mathbb{N}$, $a_{j}$ verifies:
\begin{enumerate}[(i)]
\item $a_{j}$ is uniformly bounded on a neighbourhood of $\theta_0$;
 \item  the first and second derivatives of $a_{j}$ are continuous and bounded on a neighbourhood of $\theta_0$.
\end{enumerate}
and that:
\begin{equation}
 \forall \delta>0, \exists C_l, \forall j \in \mathbb{N}^* \quad \Big \vert\frac{\partial a_j}{\partial \theta_l}(\theta_0) \Big \vert \leq C_l j^{-1+\delta} .\label{majderiveeaj}
 \end{equation}
We have then the following result:
\begin{displaymath}
\E\left( \widetilde{X'_{T,k}}(1)-\widetilde{X'_k}(1)
\right)^2= \OO\left( \frac{k^{2d}}{T}\right) .
\end{displaymath}
\end{thm}

An example to which our theorem applies is the fractionally integrated processes. In this case, the parameter $\theta$ is scalar and corresponds to the long-memory parameter $d$. Assumptions A0-A6 hold for fractionally processes. We define $d_0$ by $d_0:= \theta_0$ and then we have $a_{j}$:
\begin{displaymath}
a_{j}(d):=\left(\frac{\Gamma(j-d)}{\Gamma(j+1)\Gamma(-d)}\right)
\end{displaymath}
Since the gamma function $\Gamma$ is analytic on $\{\mathbb{C} \setminus \mathbb{N}\}$, there exists a neighbourhood of $d_0$ on which the function $a_{j}$ and its first and second derivatives are bounded.
Finally when $j \rightarrow + \infty $,
\begin{displaymath}
 \left\vert\frac{\partial a_j}{\partial d}(d_0) \right\vert \sim C j^{-1} 
\end{displaymath}
where $C$ is a constant.
As a consequence our Theroem can be applied on the class of fractionally integrated noise because all the assumptions hold. Similarly we can also show that the class of FARIMA time series verify the assumptions of this Theorem.

\begin{proof}
We first define the following vector:
\begin{displaymath}
\alpha_k^*:=\left( a_1(\widehat{\theta}_T)-a_1(\theta_0), \ldots,a_k(\widehat{\theta}_T)-a_k(\theta_0) \right) 
\end{displaymath}
where $v^*$ is the transpose vector of $v$ and 
\begin{displaymath}
\left( {\bf X}_{1}^k \right)^*:=\left( X_k, \ldots , X_{1} \right).
\end{displaymath}

\begin{eqnarray*}
\E\left( \widetilde{X'_{T,k}}(1)-\widetilde{X'_k}(1)
\right)^2&=&
\E\left( 
\left( a_1(\widehat{\theta}_T)-a_1, \ldots,a_k(\widehat{\theta}_T)-a_k \right) 
\left( 
\begin{array}{c}
X_k \\
\vdots \\
X_{1}
\end{array}
\right) 
\right)^2\\
&=& \E\left[ \left(\alpha_k^* {\bf X}_{1}^k \right)^2 \right]\\
&=& \mathrm{trace} \left( \E\left( 
 \alpha_k^*{\bf X}_{1}^k
\right)^2\right)\\
&=& \E \left(\mathrm{trace} \left( \alpha_k\alpha_k^*
{\bf X}_{1}^k\left( {\bf X}_{1}^k\right)^*
\right)  \right) \\
&=& \mathrm{trace}\left( \E\left( \alpha_k\alpha_k^*
 \right)  \Sigma_k\right)
\end{eqnarray*}
with 
\begin{displaymath}
 \Sigma_k:=\E\left({\bf X}_{1}^k\left( {\bf X}_{1}^k\right)^*\right)
\end{displaymath}

Let us first study the covariance matrix of the estimated coefficients $\E(\alpha_k\alpha_k^*)$.
We can write $\left( \alpha_k\alpha_k^*\right)_{i,j}=\E\left(g_{i,j}(\hat{\theta}) \right) $ when $g_{i,j}$ is defined by $ g_{i,j}:\theta \mapsto (a_i(\theta)-a_i(\theta_0))(a_j(\theta)-a_j(\theta_0))$.
We then use an order 2 Taylor series expansion of $g_{i,j}$ and apply Theorem 5.4.3 form \cite{dl}.We will refer to the following version. \\
If the following assumptions hold 
 \begin{enumerate}[(i)]
 \item  $\forall m \in \llbracket 1,p \rrbracket, \quad\E\left( |\widehat{\theta}_{T,m}-\theta_{0,m}|^3\right) = \OO(\eta(T))$ where $\widehat{\theta}_{T,m}$ is the $m^{th}$ entry of $\widehat{\theta}_{T}$;
 \item  $ \widehat{\theta}_T \rightarrow \theta_0, \quad \mathbb{P}\mathrm{-a.s.}$;
 \item  $g_{i,j}$ is uniformly bounded on a neighbourhood of $\theta_0$;
 \item  the first and the second derivatives $g_{i,j}$ are continuous and bounded on a neighbourhood of $\theta_0$
 \end{enumerate}
 then 
 \begin{eqnarray*}
 \E\left( g_{i,j}(\widehat{\theta}_T )\right) &=&g_{i,j}(\theta_0)+\sum_{l=1}^p \E(\widehat{\theta}_{T,l}-\theta_{0,l})\frac{\partial g_{i,j}}{\partial \theta_l}\left( \theta_{0}\right) \\
 &&+\frac{1}{2}\sum_{l=1}^p\sum_{n=1}^p\frac{\partial^2 g_{i,j} }{\partial \theta_l \partial \theta_n}(\theta_0)\E\left( \left( \widehat{\theta}_{T,l}-\theta_{0,l}\right) \left(\widehat{\theta}_{T,n}-\theta_{0,n} \right) \right) +\OO(\eta(T)).
 \end{eqnarray*}

 By assumption, conditions (ii) et (iv) hold. We note also that:
 \begin{displaymath}
 g_{i,j}(\theta_0)=0 \quad \textrm{et} \quad \forall l \in \llbracket1,p\rrbracket, \frac{\partial g_{i,j}}{\partial \theta_l}\left( \theta_{0}\right)=0
 \end{displaymath}

Next we compute the fourth order moments of $\widehat{\theta}_T-\theta_0$ in order to estimate the second and the third moments. We define:
 \begin{eqnarray*}
 \sigma_T(\theta)&:=&\left[ \frac{1}{2 \pi} \int_{-\pi}^{\pi}\left[ f(\lambda,\theta)\right]^{-1} I_T(\lambda) \dd \lambda \right]\\
 &=&\frac{{\bf{Y}}'A_T(\theta) {\bf{Y}} }{T}
 \end{eqnarray*}
where $\left( {\bf{Y}}\right)^* =(Y_1, \ldots,Y_T)$ and
 \begin{displaymath}
 \left( A_T(\theta)\right)_{j,l}:=\frac{1}{(2 \pi) ^2}\int_{-\pi}^{\pi} \mathrm{e}^{i(j-l)\lambda}\left[f(\lambda,\theta) \right]^{-1} \dd \lambda .
 \end{displaymath}

 We follows now the proof of \cite{foxtaqqu}. Since $\widehat{\theta}_T=\underset{\theta}{ \mathrm{argmin}}\left\lbrace \sigma_T(\theta)\right\rbrace $ and according to the mean-value theorem, we have:
 \begin{displaymath}
 \exists \theta^* \textrm{ such that } \vert \theta^*-\theta_0\vert  \leq  \vert\hat{\theta}-\theta_0 \vert \textrm{ and } \quad\hat{\theta}-\theta_0 = -\left[\left(  \frac{\partial^2}{\partial \theta_i\theta_j}\sigma_T(\theta^*)\right)_{1\leq i,j\leq p} \right]^{-1} \frac{\partial }{\partial \theta}\sigma_T(\theta_0).
 \end{displaymath}
 It is justified because $\theta \mapsto [f(\lambda,\theta)]^{-1}$ is twice differentiable with respect to $\theta$ and all the partial derivatives are integrable on $[-\pi,\pi]$ with respect to $\lambda$ by assumption A3. It follows from \cite{foxtaqqu} that:
 \begin{equation}
 \frac{\partial^2}{\partial \theta_i \theta_j}\sigma_T(\theta^*) \overset{\mathbb{P}\mathrm{-a.s.}}{\longrightarrow} \frac{\sigma_{\e}^2}{2 \pi}\int_{-\pi}^{\pi} \left( \frac{\partial}{\partial\theta_i}f^{-1}(\lambda,\theta_0)\right) \left( \frac{\partial}{\partial\theta_j}f^{-1}(\lambda,\theta_0)\right)f ^2(\lambda,\theta_0) \dd \lambda:= W_{i,j}
\label{defW} 
\end{equation}
 where $W:=\left(W_{i,j} \right)_{1\leq i,j\leq p} $ is a positive definite matrix. Since the matrix norm $x \mapsto \Vert x \Vert_4$ is continuous, there exists $C>0$ such that $\Vert W \Vert_4>C$ and:
 \begin{eqnarray}
\exists M \in \mathbb{N}, \quad T>M\quad, \Big\Vert\widehat{\theta}_k-\theta_0\Big\Vert_4 & \leq & C\Bigg\Vert\frac{\partial }{\partial \theta}\sigma_T(\theta_0)\Bigg\Vert_4\quad \mathbb{P}\mathrm{-a.s.}\label{majnorm4}
 \end{eqnarray}
Using this inequality, we can now estimate the fourth moments for any $m \in \llbracket 1,p \rrbracket$:
\begin{eqnarray*}
\E\left[ \Bigg\vert\frac{\partial }{\partial \theta_m}\sigma_T(\theta_0)\Bigg\vert^4\right]&= &\E\left[ \left( \frac{{\bf{Y}}'\frac{\partial A_T(\theta_0)}{\partial \theta_m} {\bf{Y}}}{T}\right) ^4\right].
 \end{eqnarray*}
 Let $m \in \llbracket 1,p \rrbracket$. We define the matrix $\Delta_m$ with $(j,l)$-th entries: $$\delta_{j,l}:= \int_{-\pi}^{\pi}\mathrm{e}^{i(j-l)\lambda}\frac{\partial }{\partial \theta_m} f ^{-1}(\lambda,\theta_0)\dd \lambda.$$

 Next we rewrite this expression as:
\begin{eqnarray}
\E\left[ \left( \frac{\partial }{\partial \theta_m}\sigma_T(\theta_0)\right)^4\right]   & = & T^{-4}\E\left(  \left(\sum_{j=1}^T\sum_{l=1}^TY_j Y_l\delta_{j,l}\right) ^4\right) \nonumber \\
&=&T^{-4}\sum_{j_1,j_2,\ldots,j_8=1}^T\delta_{j_1,j_2}\delta_{j_3,j_4}\delta_{j_5,j_6}\delta_{j_7,j_8}\E\left(Y_{j_1} Y_{j_2}\ldots Y_{j_8}\right) \label{decmoments}
\end{eqnarray}
The process is Gaussian then all the moments are a function of the autocovariances (see \cite{momentsgaussiens}). In equation \eqref{decmoments}, we can rewrite each fourth moment in the sum as a linear combination of product of 4 covariances. We then count how many covariances belongs to the set $S=\left\lbrace \E\left(Y_{j_1}Y_{j_2} \right), \E\left(Y_{j_3}Y_{j_4} \right), \E\left(Y_{j_5}Y_{j_6} \right), \E\left(Y_{j_7}Y_{j_8} \right) \right\rbrace $:
\begin{enumerate}
\item either we have $\E(Y_{j_1}Y_{j_2})\times C$ and we can distinguish the following possibilities: 
\begin{itemize}
\item $C=\E(Y_{j_3}Y_{j_4})\E\left(Y_{j_5}Y_{j_6} \right)\E\left(Y_{j_7}Y_{j_8} \right)$ only one possibility or;
\item $C$ has one element in $S$ and no other which makes 6 possibilities =(3 choices in $S$)$\times$(2 choices for the other covariances) or;
\item $C$ has no elements in $S$ which makes 8 possibilities. First choose a complement for $Y_{j_3}$ (4 possibilities) then a complement for $Y_{j_4}$ (only 2 possibilities because the pairs in $S$ are excluded);
\end{itemize}
\item or $Y_{j_1}$ is with $Y_{j_l}$, $l>2$, which makes 5 possibilities. Let us assume that $Y_{j_1}$ is associated with $Y_{j_3}$. We can then distinguish the following cases:
\begin{itemize}
\item we obtain 2 pairs in $S$ which are consequently $\left(Y_{j_5},Y_{j_6} \right) $ and $\left(Y_{j_7}, Y_{j_8}\right) $ or;
\item we have only one couple in $S$ which makes (2 choices in $S$) $\times ((\mathrm{C}^2_4-1)$ choices for the other covariances) or;
\item we have non elements in $S$: either $Y_{j_2}$ is the complement of $Y_{j_4}$ and then we have only 2 possibilities, or we have 4 choices for the complement of $Y_{j_2}$ and only 2 for $Y_{j_4}$. Finally we have 10 possibilities.
\end{itemize}
\end{enumerate}

Therefore we obtain:
\begin{eqnarray*}
&&\E\left[ \left( \frac{\partial }{\partial \theta_m}\sigma_T(\theta_0)\right)^4\right] \nonumber \\  & = & T^{-4}\left( \sum_{j=1}^T\sum_{l=1}^T \delta_{j,l} \sigma(j-l)\right) ^4  \nonumber \\
&+&T^{-4}6 \left( \sum_{j=1}^T\sum_{l=1}^T \delta_{j,l} \sigma(j-l)\right) ^2\sum_{j_1,j_2,j_3,j_4=1}^T \delta_{j_1,j_2} \delta_{j_3,j_4} \sigma(j_1-j_3)\sigma(j_2-j_4)\nonumber \\
&+&T^{-4}8\left( \sum_{j=1}^T\sum_{l=1}^T \delta_{j,l} \sigma(j-l)\right)\sum_{j_1,j_2,j_3,j_4,j_5,j_6=1}^T \delta_{j_1,j_2} \delta_{j_3,j_4}  \delta_{j_5,j_6}\sigma(j_1-j_3)\sigma(j_4-j_5)\sigma(j_6-j_2)\nonumber \\
&+&5T^{-4}\left( \sum_{j=1}^T\sum_{l=1}^T \delta_{j,l} \sigma(j-l)\right) ^2\sum_{j_1,j_2,j_3,j_4=1}^T \delta_{j_1,j_2} \delta_{j_3,j_4} \sigma(j_1-j_3)\sigma(j_2-j_4)\nonumber \\
&+&5T^{-4}10\left( \sum_{j=1}^T\sum_{l=1}^T \delta_{j,l} \sigma(j-l)\right)\sum_{j_1,j_2,j_3,j_4,j_5,j_6=1}^T \delta_{j_1,j_2} \delta_{j_3,j_4}  \delta_{j_5,j_6}\sigma(j_1-j_3)\sigma(j_4-j_5)\sigma(j_6-j_2)\nonumber \\
&+&5T^{-4}10\sum_{j_1,j_2,j_3,j_4,j_5,j_6,j_7,j_8=1}^T \delta_{j_1,j_2} \delta_{j_3,j_4}  \delta_{j_5,j_6}\delta_{j_7,j_8}\sigma(j_1-j_3)\sigma(j_4-j_5)\sigma(j_6-j_7)\sigma(j_8-j_2)\nonumber
\end{eqnarray*}

All the terms of this sum are like:
\begin{equation}
\sum_{j_1,\ldots,j_{2p}=1}^T\delta_{j_1,j_2}\ldots \delta_{j_{2p-1},j_{2p}}\sigma(j_1-j_3)\ldots\sigma(j_{2p}-j_2):=S_{p,T}. \label{general}
\end{equation}
Note that $S_{p,T}=\mathrm{trace}\left( \left(\Sigma_T \Delta_m \right)^p \right)$ and that $\Delta_m$ is the covariance matrix defined by the spectral density:
\begin{displaymath}
\frac{\partial }{\partial \theta_m} f^{-1}(\lambda,\theta_0)
=\OO\left( \lambda^{\alpha(\theta_0)-\delta}\right) \quad \textrm{as} \quad\lambda \rightarrow 0,\quad \textrm{for any } \delta>0
\end{displaymath}
by assumption A3.
By applying the Theorem 1 of \cite{centrallimit}, we prove that:
\begin{eqnarray}
&&\lim_{T \rightarrow + \infty}\frac{1}{T}\sum_{j_1,\ldots,j_{2p}=1}^T \delta_{j_1,j_2}\ldots \delta_{j_{2p-1},j_{2p}}\sigma(j_2-j_3)\ldots\sigma(j_{2p}-j_1) \nonumber\\
&=& \left( 2  \pi \right)^{2p-1} \int_{-\pi}^{\pi}\left(\frac{\partial }{\partial \theta_m} f^{-1}(\lambda,\theta_0)  f(\lambda,\theta_0)\right)^p  \dd \lambda\nonumber\\
&=&\OO(1).\label{maj1} 
\end{eqnarray}
It follows from assumptions A2 and A3 that this integral is always finite. We need a more precise result for the term:
\begin{displaymath}
 \sum_{j=1}^T\sum_{l=1}^T \delta_{j,l} \sigma(j-l).
 \end{displaymath}
 Despite this can also be expressed like \eqref{general}, the estimate given below is not sufficient to conclude.
By \cite{foxtaqqu}[proof of Theorem 2], we have:
\begin{equation}
\forall \delta>0, \quad \sum_{j=1}^T\sum_{l=1}^T \delta_{j,l} \sigma(j-l) = \OO\left( T^{\delta}\right).\label{maj2}
\end{equation}
By \eqref{maj1} et \eqref{maj2}, we may conclude that 
\begin{equation}
\forall \delta>0, \quad \E\left[ \left( \frac{\partial }{\partial \theta_m}\sigma_T(\theta_0)\right)^4\right]=\OO(T^{-3+\delta}).\label{moment4dsigma}
\end{equation} 
Next using the asymptotic estimate of the fourth moments, we can now obtain asymptotic properties for the second moments:
\begin{displaymath}
\E\left[\left( \widehat{\theta}_{T,j}-\theta_{0,j}\right) \left(\widehat{\theta}_{T,l}-\theta_{0,l} \right)  \right]. 
\end{displaymath}
First we have to prove the uniform integrability of $ \sqrt{T}\left(\widehat{\theta}_{T,j}-\theta_{0,j}\right) \sqrt{T}\left(\widehat{\theta}_{T,l}-\theta_{0,l} \right) $:
\begin{eqnarray*}
T^2\E\left(\left( \widehat{\theta}_{T,j}-\theta_{0,j}\right)^2 \left(\widehat{\theta}_{T,l}-\theta_{0,l} \right)^2  \right)& \leq & T^2\sqrt{\E\left(\left( \widehat{\theta}_{T,j}-\theta_{0,j}\right)^4\right)\E\left(\left( \widehat{\theta}_{T,l}-\theta_{0,l}\right)^4\right)}\\
& \leq &T^2\E\left(\Big\Vert\widehat{\theta}_{T} -\theta_{0}\Big\Vert_4^4\right) \\
& \leq &T^2C^4\E\left(\Bigg\Vert\frac{\partial }{\partial \theta}\sigma_T(\theta_0)\Bigg\Vert_4^4\right)
\end{eqnarray*}
if $T>M$ from \eqref{majnorm4}. 
By applying result \eqref{moment4dsigma}, we conclude that:
\begin{eqnarray*}
\forall \delta>0, \quad T^2\E\left(\left( \widehat{\theta}_{T,j}-\theta_{0,j}\right)^2 \left(\widehat{\theta}_{T,l}-\theta_{0,l} \right)^2  \right)& = &\OO\left( T^{2-3+\delta}\right) \\
 &= &\OO\left( 1\right)
\end{eqnarray*}
We have proved the uniform integrability of $ \sqrt{T}\left(\widehat{\theta}_{T,j}-\theta_{0,j}\right) \sqrt{T}\left(\widehat{\theta}_{T,l}-\theta_{0,l} \right) $ since if $\E\left(X_T^2 \right)$ is finite for any $T$,then the collection $(X_T)$ is uniformly integrable. Moreover according to \cite{foxtaqqu} [Theorem 2]:
\begin{displaymath}
\sqrt{T}\left( \widehat{\theta}_{T} -\theta_{0}\right) \xrightarrow{\mathcal{L}} \mathcal{N}\left(0,4 \pi W^{-1} \right) 
\end{displaymath}
where $W$ is the matrix defined in \eqref{defW} and we have also the following convergence in law:
\begin{displaymath}
h_{j,l}\left(\sqrt{T}\left( \widehat{\theta}_{T} -\theta_{0}\right)  \right):=T\left(\widehat{\theta}_{T,j}-\theta_{0,j}\right) \left(\widehat{\theta}_{T,l}-\theta_{0,l} \right) \xrightarrow{\mathcal{L}} h_{j,l}\left(Z\right) 
\end{displaymath}
with $Z \sim  \mathcal{N}\left(0,4 \pi W^{-1} \right)$.
By the convergence in law and the uniform integrability we apply Theorem 5.4 in \cite{billingsley} that:
 \begin{displaymath}
 \E\left( T\left(\widehat{\theta}_{T,j}-\theta_{0,j}\right) \left(\widehat{\theta}_{T,l}-\theta_{0,l} \right)\right) \rightarrow 4 \pi W^{-1}_{j,l} \textrm{ as } T \rightarrow + \infty.
 \end{displaymath}

 Now we give an asymptotic bound for the third order moment by applying the Cauchy-Schwarz inequality. Using the inequalities \eqref{majnorm4} and \eqref{moment4dsigma}, we conclude that:
 \begin{eqnarray*}
 \exists C>0, \:\E\left(\left\vert \widehat{\theta}_{T,j}-\theta_{0,j}\right\vert^3 \right) & \leq & \sqrt{\E\left(\widehat{\theta}_{T,j}-\theta_{0,j}\right)^2\E\left(\widehat{\theta}_{T,j}-\theta_{0,j}\right)^4} \\
 & \leq &\sqrt{C T^{-1} T^{-3+\delta} }, \quad \forall  \delta>0\\
 &=&\OO \left( T^{-2+\delta}\right) , \quad \forall  \delta>0
 \end{eqnarray*}
 We obtain the following Taylor series for any $\delta>0 $:
 \begin{eqnarray*}
 \E\left( (a_i(\widehat{\theta}_T)-a_i(\theta_0))(a_j(\widehat{\theta}_T)-a_j(\theta_0))\right) &=&2\pi T^{-1}\sum_{l=1}^m \sum_{n=1}^m \frac{\partial^2 g_{i,j} }{\partial \theta_l \partial \theta_n}(\theta_0)W^{-1}_{l,n}+\OO \left( T^{-2+\delta}\right) \\
 &\sim&2\pi T^{-1}\sum_{l=1}^m \sum_{n=1}^m\left[ \frac{\partial a_i(\theta_0)}{\partial \theta_l}\frac{\partial a_j(\theta_0)}{\partial \theta_n}+\frac{\partial a_j(\theta_0)}{\partial \theta_l}\frac{\partial a_i(\theta_0)}{\partial \theta_n}\right] W^{-1}_{l,n}.
 \end{eqnarray*}
 
 We can now conclude and find an asymptotic equivalent of $\E\left( \alpha_k\alpha_k^*\right) $. Since $W^{-1}$ is symmetric:
 \begin{eqnarray}
 \E\left( \alpha_k\alpha_k^*\right) \sim 4\pi T^{-1} D W^{-1}  D^* \label{alphak}
 \end{eqnarray}
with
\begin{displaymath}
D:= \left( 
\begin{array}{ccc}
 \frac{\partial a_1(\theta_0)}{\partial \theta_1} &\ldots & \frac{\partial a_1(\theta_0)}{\partial \theta_m}\\
 \vdots & \ddots & \vdots \\
 \frac{\partial a_k(\theta_0)}{\partial \theta_1} &\ldots & \frac{\partial a_k(\theta_0)}{\partial \theta_m}
\end{array}
\right) .
\end{displaymath}
 $W^{-1}$ is a positive definite matrix because $W$ is too. So it can be expressed as :
\begin{displaymath}
W^{-1}=P^* \left(\begin{array}{cccc}
\lambda_1 &0&\ldots& 0\\
0& \lambda_2 & \ddots &\vdots \\
\vdots &\ddots & \ddots &\vdots \\
0 & \ldots &  0& \lambda_m
\end{array}
 \right) P
\end{displaymath}
 where $P=(p_{ij})_{1 \leq i,j \leq m}$ is an orthogonal matrix and the $(\lambda_i)_{1 \leq i \leq m}$ are the positive eigenvalues of $W^{-1}$. We may rewrite our expression as:
 \begin{eqnarray}
  D W^{-1}  D^*&=& \sum_{r=1}^m \left( \left[\sqrt{\lambda_r} \sum_{l=1}^m p_{rl} \frac{\partial a_i}{\partial \theta_l}(\theta_0) \right] \left[\sqrt{\lambda_r} \sum_{l=1}^m p_{rl} \frac{\partial a_j}{\partial \theta_l}(\theta_0) \right] \right)_{1 \leq i,j \leq k} \nonumber\\
  &=&\sum_{r=1}^m  \beta_r^* \beta_r \label{betar}
 \end{eqnarray}
where $\beta_r$ is the vector $\left( \sqrt{\lambda_r} \sum_{l=1}^m p_{rl} \frac{\partial a_i}{\partial \theta_l}(\theta_0)\right)_ {1 \leq i \leq k}$.

\begin{eqnarray*}
 \E\left( 
\left( a_1(\widehat{\theta}_T)-a_1(\theta_0), \ldots,a_k(\widehat{\theta}_T)-a_k (\theta_0)\right) 
\left( 
\begin{array}{c}
X_0 \\
\vdots \\
X_{-k+1}
\end{array}
\right) 
\right)^2 &=& \mathrm{trace}\left( \E\left( \alpha_k\alpha_k^*
 \right)  \Sigma_k\right)\\
 &\sim& 4 \pi T^{-1}\sum_{r=1}^m \mathrm{trace}\left( \beta_r^* \beta_r\Sigma_k\right)
 \end{eqnarray*}
from \eqref{alphak} and \eqref{betar}. Therefore we obtain:
 \begin{eqnarray*}
 \mathrm{trace}\left( \E\left( \alpha_k\alpha_k^*\right)  \Sigma_k\right)&\sim& 4 \pi T^{-1}\sum_{r=1}^m \beta_r\Sigma_k \beta_r^* \\
 &\leq&  4 \pi T^{-1}\sum_{r=1}^m \Lambda_k \Vert\beta_r\Vert^2_2
 \end{eqnarray*}
where $\Lambda_k$ is the greatest eigenvalue of $\Sigma_k$. The last inequality is a consequence of $\Sigma_k$ being symmetric matrix.
Following assumption \eqref{majderiveeaj}, we have:
 \begin{displaymath}
 \forall \delta>0, \exists C_l, \forall j \in \mathbb{N}^* \quad \Big\vert\frac{\partial a_j}{\partial \theta_l}(\theta_0) \Big\vert\leq C_l j^{-1+\delta}
 \end{displaymath}
and we can hence estimate $\Vert\beta_r\Vert^2_2$. Let $\delta= 1/2$, there exists $C_1, \ldots,C_m$ such that:
\begin{eqnarray*}
\Vert\beta_r\Vert^2_2&=& \sum_{j=1}^k \lambda_r \sum_{l_1=1}^m\sum_{l_2=1}^m p_{rl_1}p_{rl_2}\frac{\partial a_j}{\partial \theta_{l_1}}(\theta_0)\frac{\partial a_j}{\partial \theta_{l_2}}(\theta_0)\\
& \leq & \lambda_r \sum_{l_1=1}^m\sum_{l_2=1}^m \left|p_{rl_1}p_{rl_2}\right|\sum_{j=1}^k C_{l_1} C_{l_2}j^{-3}\\
& \leq & \lambda_r \sum_{l_1=1}^m\sum_{l_2=1}^m\left|p_{rl_1}p_{rl_2}\right|\sum_{j=1}^{+ \infty}C_{l_1} C_{l_2}j^{-3}:=C_r(\theta_0)
\end{eqnarray*}
where $C_r(\theta_0)$ does not depend on $k$.

From \cite{boettcher}, the spectral norm of a Toeplitz matrix (its spectral norm), whose symbol has the form $\lambda \mapsto \lambda^{-\alpha}L(\lambda)$ with $L$ is a bounded, continuous at 0 function and does not vanish at $0$, is equivalent to $Ck^{\alpha}$ with $C$ constant. We conclude the proof:
\begin{displaymath}
\E\left( 
\left( a_1(\widehat{\theta}_T)-a_1(\theta_0), \ldots,a_k(\widehat{\theta}_T)-a_k(\theta_0) \right) 
\left( 
\begin{array}{c}
X_0 \\
\vdots \\
X_{-k+1}
\end{array}
\right) 
\right)^2 \leq  C 4\pi \sum_{r=1}^{m}C_r(\theta_0)  \frac{k^{\alpha(\theta_0)}}{T}
\end{displaymath}
with $C$ constant.
\end{proof}

\subsection{Conclusion}
\label{conclusion1}
\par Prediction with the Wiener-Kolmogorov predictor involves two mean-squared error components: the first is due to the truncation to $k$ terms and this is bounded by $\OO(k^{-1})$, the second is due to the estimation of the coefficients $a_j$ from a realisation of the process of length $T$ and is bounded by $\OO(k^{2d}/T)$. The mean-squared difference between the best linear predictor $\widetilde{X_k}(1)$ and our predictor is given by:
\begin{eqnarray*}
\Big \Vert \widetilde{X_{T,k}}(1)-\widetilde{X_k}(1) \Big \Vert_2
 &\leq &\Big \Vert \widetilde{X_{T,k}}(1)-\widetilde{X_k'}(1) \Big \Vert_2+\Big \Vert \widetilde{X_k'}(1)-\widetilde{X_k}(1) \Big \Vert_2 \\
&\leq&\OO(k^{-1/2})+\OO(k^{d}/\sqrt{T}).
\end{eqnarray*}
If we want to compare the two types of prediction errors, we need a relation between the rate of convergence of $T$ and $k$ to $+ \infty$. For example, if $T=\mathrm{o}(k^{2d+1})$, the error due to the estimation of the coefficients is predominant and gives the bound for the general error.

\paragraph*{} Truncating to $k$ terms the series which defines the Wiener-Kolmogorov predictor amounts to using an AR($k$) model for predicting. Therefore in the following section we look for the AR($k$) which minimizes the forecast error.

\section{The Autoregressive Models Fitting Approach} \label{sectionarp}

In this section we shall develop a generalisation of the "autoregressive model fitting" approach developed by \cite{arp} in the case of fractionally integrated noise F($d$) (defined in \eqref{Fd}). We study asymptotic properties of the forecast mean-squared error when we fit a misspecified AR($k$) model to the long-memory time series $(X_n)_{n \in \mathbb{Z}}$. 

\subsection{Rationale}

Let $\Phi$ a $k^{\textrm{th}}$ degree polynomial defined by: 
\begin{displaymath}
\Phi(z)=1-a_{1,k} z-\ldots -a_{k,k} z^k.
\end{displaymath}
We assume that $\Phi$ has no zeroes on the unit disk. We define the process $(\eta_n)_{n\in \mathbb{Z}}$ by:
\begin{displaymath}
\forall n \in \mathbb{Z}\textrm{, }\eta_n =\Phi(B) X_n
\end{displaymath} where $B$ is the backward shift operator.
Note that $(\eta_n)_{n\in \mathbb{Z}}$ is not a white noise series because $(X_n)_{n\in \mathbb{Z}}$ is a long-memory process and hence does not belong to the class of autoregressive processes. Since $\Phi$ has no root on the unit disk, $(X_n)_{n\in \mathbb{Z}}$ admits a moving-average representation as the fitted AR($k$) model in terms of $(\eta_n)_{n\in \mathbb{Z}}$:
\begin{displaymath}
X_n=\sum_{j=0}^{\infty} c(j)\eta_{n-j}.
\end{displaymath}
If $(X_n)_{n\in \mathbb{Z}}$ was an AR($k$) associated with the polynomial $\Phi$, the best next step linear predictor would be:
\begin{eqnarray*}
\widehat{X_n}(1)&=& \sum_{j=1}^{\infty} c(i) \eta_{t+1-i}\\
&=& a_{1,k}X_n +\ldots+a_{k,k} X_{n+1-k} \textrm{ si }n \geqslant k.
\end{eqnarray*}
Here $(X_n)_{n\in \mathbb{Z}}$ is a long-memory process which verifies the assumptions of Section \ref{intro}. Our goal is to express the polynomial $\Phi$ which minimizes the forecast error and to estimate this error.

\subsection{Mean-Squared Error}
\label{yaj}
There exists two approaches in order to define the coefficients of the $k^{\mathrm{th}}$ degree polynomial $\Phi$: the spectral approach and the time approach. 
\par In the time approach, we choose to define the predictor as the projection mapping on to the closed span of the subset $\{X_n,\ldots,X_{n+1-k}\}$ of the Hilbert space L$^2(\Omega,\mathcal{F},\mathbb{P})$ with inner product $<X,Y>=\E(XY)$. Consequently the coefficients of $\Phi$ verify the equations, which are called the $k^{\mathrm{th}}$ order Yule-Walker equations:
\begin{equation}
\forall j \in \llbracket1,k \rrbracket, \quad \sum_{i=1}^k a_{i,k}\sigma (i-j)=\sigma(j) \label{yw}
\end{equation}

\par  The mean-squared prediction error is:
\begin{displaymath}
\mathbb{E}\big[ \big(\widehat{X_n}(1)-X_{n+1}\big)²\big]=c(0)²\mathbb{E}(\eta_{n+1}^2)=\mathbb{E}(\eta_{n+1}^2).
\end{displaymath}

We may write the moving average representation of $(\eta_n)_{n\in \mathbb{N}}$ in terms of $(\varepsilon_n)_{n\in \mathbb{N}}$:
\begin{eqnarray*}
\eta_n&=& \sum_{j=0}^{\infty} \sum_{k=0}^{min(j,p)}\Phi_k b(j-k) \varepsilon_{n-j}\\
&=& \sum_{j=0}^{\infty} t(j)\varepsilon_{n-j}
\end{eqnarray*} 
with
\begin{displaymath}
\forall j \in \mathbb{N}, \quad t(j)=\sum_{k=0}^{min(j,p)}\Phi_k b(j-k).
\end{displaymath}

Finally we obtain:
\begin{displaymath}
\mathbb{E}\big[ \big(\widehat{X_n}(1)-X_{n+1}\big)²\big]=\sum_{j=0}^{\infty} t(j)² \sigma^2_{\varepsilon}.
\end{displaymath}
In the spectral approach, minimizing the prediction error is equivalent to minimizing a contrast between two spectral densities:
\begin{equation*}
\int_{-\pi}^{\pi}\frac{f(\lambda)}{g(\lambda,\Phi)}\dd \lambda \label{distance}
\end{equation*} 
where $f$ is the spectral density of $X_n$ and $g(.,\Phi)$ is the spectral density of the AR(p) process defined by the polynomial $\Phi$ (see for example \cite{dist}),so:
\begin{align*}
\int_{-\pi}^{\pi}\frac{f(\lambda)}{g(\lambda,\Phi)}\dd \lambda &= \int_{-\pi}^{\pi}\Big|\sum_{j=0}^{\infty}b(j)e^{-ij\lambda}\Big|²\Big|\Phi(e^{-i\lambda })\Big|²\dd \lambda\\
&=\int_{-\pi}^{\pi}|\sum_{j=0}^{\infty}t(j)e^{-ij\lambda}|²\dd \lambda\\
&=2\pi\sum_{j=0}^{\infty}t(j)².
\end{align*}
In both approaches we nedd to minimize $\sum_{j=0}^{\infty}t(j)$.

\subsection{Rate of Convergence of the Error by AR($k$) Model Fitting}

In the next theorem we derive an asymptotic expression for the prediction error by fitting autoregressive models to the series:
\begin{thm}
Assume that $(X_n)_{n \in \mathbb{Z}}$ is a long-memory process which verifies the assumptions of Section \ref{intro}. If $0<d<\frac{1}{2}$:
\begin{displaymath}
\mathbb{E}\big[ \big(\widehat{X_k}(1)-X_{k+1}\big)²\big]-\sigma_{\e}^2=\OO(k^{-1})
\end{displaymath}
\end{thm}
\begin{proof}
Since fitting an AR($k$) model minimizes the forecast error using $k$ observations, the error by using truncation way is bigger. So, since the truncation method involves an error bounded by $\OO\left(k^{-1} \right) $, we obtain:
\begin{displaymath}
\mathbb{E}\big[ \big(\widehat{X_k}(1)-X_{k+1}\big)²\big]-\sigma_{\e}^2=\OO(k^{-1}).
\end{displaymath}
Consequently we only need to prove that this rate of convergence is attained
. This is the case for the fractionally integrated processes defined in \eqref{Fd}. We want the error made when fitting an AR($k$) model in terms of the Wiener-Kolmogorov truncation error. Note first that the variance of the white noise series is equal to:
\begin{displaymath}
\sigma_{\e}^2= \int_{-\pi}^{\pi}f(\lambda) \left| \sum_{j=0}^{+\infty}a_j e^{ij\lambda} \right|^2 \dd \lambda.
\end{displaymath}
Therefore in the case of a fractionally integrated process F($d$) we need only show that:
\begin{displaymath}
\int_{-\pi}^{\pi}f(\lambda) \left| \sum_{j=0}^{+\infty}a_j e^{ij\lambda} \right|^2 \dd \lambda-\frac{\sigma_{\e}^2}{2\pi}\int_{-\pi}^{\pi}\frac{f(\lambda)}{g(\lambda,\Phi_k)} \dd \lambda\sim C(k^{-1}).
\end{displaymath}
\begin{eqnarray*}
\int_{-\pi}^{\pi}f(\lambda) \left| \sum_{j=0}^{+\infty}a_j e^{ij\lambda} \right|^2 \dd \lambda-\frac{\sigma_{\e}^2}{2\pi}\int_{-\pi}^{\pi}\frac{f(\lambda)}{g(\lambda,\Phi_k)}\dd \lambda&=& \int_{-\pi}^{\pi}f(\lambda) \left( \left| \sum_{j=0}^{+\infty}a_j e^{ij\lambda} \right|^2- \left| \sum_{j=0}^{k}a_{j,k} e^{ij\lambda} \right|^2\right) \dd \lambda \\
&=& \sum_{j=0}^{+\infty}\sum_{l=0}^{+\infty}\left(a_ja_l-a_{j,k}a_{l,k}\right) \sigma(j-l)
\end{eqnarray*}
we set $a_{j,k}=0$ if $j>k$.
\begin{eqnarray}
&&\sum_{j=0}^{+\infty}\sum_{l=0}^{+\infty}\left(a_j a_l -a_{j,k}a_{l,k}\right) \sigma(j-l) \\ 
&=&\sum_{j=0}^{+\infty}\sum_{l=0}^{+\infty}(a_ja_l-a_{j,k}a_l)\sigma(j-l) 
+\sum_{j=0}^{+\infty}\sum_{l=0}^{+\infty}(a_{j,k}a_l-a_{j,k}a_{l,k})\sigma(j-l) \nonumber\\
&=&\sum_{j=0}^{+\infty}(a_j-a_{j,k})\sum_{l=0}^{+\infty}a_l \sigma(l-j)
+\sum_{j=0}^{k}a_{j,k}\sum_{l=0}^{+\infty}(a_l-a_{l,k})\sigma(j-l)  \quad \label{somme}
\end{eqnarray}
We first study the first term of the sum \eqref{somme}. For any  $j>0$ , we have  $\sum_{l=0}^{+\infty}a_l \sigma(l-j)=0$:
\begin{eqnarray*}
\e_n&=&\sum_{j=0}^{\infty} a_l X_{n-l}\\
X_{n-j} \e_n&=&\sum_{l=0}^{\infty} a_lX_{n-l} X_{n-j} \\
\E\left( X_{n-j} \e_n\right) &=&\sum_{l=0}^{\infty} a_l\sigma(l-j) \\
\E\left( \sum_{l=0}^{\infty} b_l \e_{n-j-l}\e_n\right)&=&\sum_{l=0}^{\infty} a_l\sigma(l-j)
\end{eqnarray*}
and we conclude that $\sum_{l=0}^{+\infty}a_l \sigma(l-j)=0$ because $(\varepsilon_n)_{n\in \mathbb{Z}}$ is an uncorrelated white noise.
We can thus rewrite the first term of \eqref{somme} like:
\begin{eqnarray}
\sum_{j=0}^{+\infty}(a_j-a_{j,k})\sum_{l=0}^{+\infty}a_l \sigma(l-j)&=&(a_0-a_{0,k})\sum_{l=0}^{+\infty}a_l \sigma(l) \nonumber \\
&=&0 \nonumber 
\end{eqnarray}
since $a_0=a_{0,k}=1$ according to definition. 
Next we study the second term of the sum \eqref{somme}:
\begin{eqnarray}
\sum_{j=0}^{k}a_{j,k}\sum_{l=0}^{+\infty}(a_l-a_{l,k})\sigma(j-l) \nonumber .
\end{eqnarray}
And we obtain that:
\begin{eqnarray}
\sum_{j=0}^{k}a_{j,k}\sum_{l=0}^{+\infty}(a_l-a_{l,k})\sigma(j-l)&=& \sum_{j=1}^{k}(a_{j,k}-a_j)\sum_{l=1}^{k}(a_l-a_{l,k})\sigma(j-l) \nonumber \\
&&+\sum_{j=1}^{k}(a_{j,k}-a_j)\sum_{l=k+1}^{+\infty}a_l\sigma(j-l) \label{milieu1}\\
&&+\sum_{j=0}^{k}a_j\sum_{l=1}^{k}(a_l-a_{l,k})\sigma(j-l) \label{milieu2}\\
&&+\sum_{j=0}^{k}a_j\sum_{l=k+1}^{+\infty}a_l\sigma(j-l) \nonumber
\end{eqnarray}
Similarly we rewrite the term \eqref{milieu1} using the Yule-Walker equations:
\begin{displaymath}
\sum_{j=1}^{k}(a_{j,k}-a_j)\sum_{l=k+1}^{+\infty}a_l\sigma(j-l)=-\sum_{j=1}^{k}(a_{j,k}-a_j)\sum_{l=0}^{k}a_l\sigma(j-l)
\end{displaymath} 
We then remark that this is equal to \eqref{milieu2}.
Hence it follows that:
\begin{eqnarray}
\sum_{j=0}^{k}a_{j,k}\sum_{l=0}^{+\infty}(a_l-a_{l,k})\sigma(j-l)&=&\sum_{j=1}^{k}(a_{j,k}-a_j)\sum_{l=1}^{k}(a_l-a_{l,k})\sigma(j-l) \nonumber \\
&&+2\sum_{j=1}^{k}(a_{j,k}-a_j) \sum_{l=k+1}^{+\infty}a_l\sigma(j-l) \nonumber \\
&&+\sum_{j=0}^{k}a_j\sum_{l=k+1}^{+\infty}a_l\sigma(j-l) \label{commeavant}
\end{eqnarray}
On a similar way we can rewrite the third term of the sum \eqref{commeavant} using Fubini Theorem:
\begin{displaymath}
\sum_{j=0}^{k}a_j\sum_{l=k+1}^{+\infty}a_l\sigma(j-l)=-\sum_{j=k+1}^{+\infty}\sum_{l=k+1}^{+\infty}a_ja_l \sigma(j-l).
\end{displaymath}
This third term is therefore equal to the forecast error in the method of prediction by truncation.

  In order to compare the prediction error by truncating the Wiener-Kolmogorov predictor and by fitting an autoregressive model to a fractionally integrated process F($d$), we need the sign of all the components of the sum \eqref{commeavant}. For a fractionally integrated noise, we know the explicit formula for $a_j$ and $\sigma(j)$:
\begin{displaymath}
\forall j>0,\quad a_j=\frac{\Gamma(j-d)}{\Gamma(j+1)\Gamma(-d)}<0 \textrm{ and } \forall j\geq 0,\quad \sigma(j)=\frac{(-1)^j\Gamma(1-2d)}{\Gamma(j-d+1)\Gamma(1-j-d)} \sigma_\e^2>0.
\end{displaymath}
In order to get the sign of $a_{j,k}-a_j$ we use the explicit formule given in \cite{coefftronccourtemem} and we easily obtain that $a_{j,k}-a_j$ is negative for all $j \in \llbracket 1, k \rrbracket$.
\begin{eqnarray*}
a_j-a_{j,k}&=&\frac{\Gamma(j-d)}{\Gamma(j+1)\Gamma(-d)}-\frac{\Gamma(k+1)\Gamma(j-d) \Gamma(k-d-j+1)}{\Gamma(k-j+1)\Gamma(j+1)\Gamma(-d)\Gamma(k-d+1)}\\
&=&-a_j\left(-1+\frac{\Gamma(k+1)\Gamma(k-d-j+1)}{\Gamma(k-j+1)\Gamma(k-d+1)} \right) \\
&=&-a_j\left(\frac{k...(k-j+1)}{(k-d)...(k-d-j+1)}-1\right) \\
&>&0
\end{eqnarray*}
since $\forall j \in \mathbb{N}^*$ $a_j<0$.
 To give an asymptotic equivalent for the prediction error, we use the sum given in \eqref{commeavant}. We have the sign of the three terms: the first is negative, the second is positive and the last is negative. Moreover the third is equal to the forecast error by truncation and we have proved that this asymptotic equivalent has order $\OO(k^{-1})$. The  prediction error by fitting an autoregressive model converges faster to 0 than the error by truncation only if the second term is equivalent to $Ck^{-1}$, with $C$ constant. Consequently, we search for a bound for $a_j-a_{j,k}$ given the explicit formula for these coefficients (see for example \cite{coefftronccourtemem}):
\begin{eqnarray*}
a_j-a_{j,k}&=&\frac{\Gamma(j-d)}{\Gamma(j+1)\Gamma(-d)}-\frac{\Gamma(k+1)\Gamma(j-d) \Gamma(k-d-j+1)}{\Gamma(k-j+1)\Gamma(j+1)\Gamma(-d)\Gamma(k-d+1)}\\
&=&-a_j\left(-1+\frac{\Gamma(k+1)\Gamma(k-d-j+1)}{\Gamma(k-j+1)\Gamma(k-d+1)} \right) \\
&=&-a_j\left(\frac{k...(k-j+1)}{(k-d)...(k-d-j+1)}-1\right) \\
&=&-a_j\left(\prod_{m=0}^{j-1} \left( \frac{1-\frac{l}{k}}{1-\frac{l+d}{k}}\right) -1\right) \\
&=&-a_j\left(\prod_{m=0}^{j-1}\left( 1+\frac{\frac{d}{k}}{1-\frac{d+l}{k}}\right) -1\right) .
\end{eqnarray*}
Then we use the following inequality:
\begin{displaymath}
\forall x \in \mathbb{R}, \quad 1+x\leq \exp(x)
\end{displaymath}
which gives us:
\begin{eqnarray*}
a_j-a_{j,k}&\leq&-a_j\left(\exp\left(\sum_{m=0}^{j-1} \frac{\frac{d}{k}}{1-\frac{d+l}{k}}\right) -1\right) \\
&\leq&-a_j\left(\exp \left(d\sum_{m=0}^{j-1}\frac{1}{k-d-l} \right) -1\right)\\
&\leq&-a_j\exp \left(d\sum_{m=0}^{j-1}\frac{1}{k-d-l} \right) 
\end{eqnarray*}
According to the previous inequality, we have:
\begin{eqnarray*}
\sum_{j=1}^{k}(a_j-a_{j,k})\sum_{l=k+1}^{+\infty}-a_l \sigma(j-l)&=&\sum_{j=1}^{k-1}(a_j-a_{j,k})\sum_{l=k+1}^{+\infty}-a_l \sigma(j-l) \\
&&+(a_k-a_{k,k})\sum_{l=k+1}^{+\infty}-a_l \sigma(k-l) \\
&\leq & \sum_{j=1}^{k-1}-a_j\exp \left(d\sum_{m=0}^{j-1}\frac{1}{k-d-m}\right)  \sum_{l=k+1}^{+\infty}-a_l \sigma(j-l) \\
&&+(-a_k)\exp \left(d\sum_{m=0}^{k-1}\frac{1}{k-d-m} \right)\sum_{l=k+1}^{+\infty}-a_l \sigma(k-l) \\
&\leq & \sum_{j=1}^{k-1}-a_j\exp \left(d\int_0^{j}\frac{1}{k-d-m} \dd m\right)\sum_{l=k+1}^{+\infty}-a_l \sigma(j-l) \\
&&+(-a_k) k^{\frac{3}{2}d}\sum_{l=k+1}^{+\infty}-a_l \sigma(k-l) 
\end{eqnarray*}
As the function $x \mapsto \frac{1}{k-d-x}$ is increasing, we use the Integral Test Theorem. The inequality on the second term follows from:
\begin{eqnarray*}
\sum_{m=0}^{k-1}\frac{1}{k-d-m} &\sim& \ln(k)\\
&\leq& \frac{3}{2}\ln(k) 
\end{eqnarray*}for $k$ large enough.
Therefore there exists $K$ such that for all $k \geq K$:
\begin{eqnarray*}
\sum_{j=1}^{k}(a_j-a_{j,k})\sum_{l=k+1}^{+\infty}-a_l \sigma(j-l)&\leq & \sum_{j=1}^{k-1}-a_j\exp \left(d \ln\left( \frac{k-d}{k-d-j}\right) \right) \sum_{l=k+1}^{+\infty}-a_l \sigma(j-l)\\
&&+(-a_k) k^{\frac{3}{2}d}\sum_{l=k+1}^{+\infty}-a_l \sigma(0)\\
&\leq &C (k-d)^d\sum_{j=1}^{k-1}j^{-d-1}(k-d-j)^{-d}\sum_{l=k+1}^{+\infty}l^{-d-1}(l-j)^{2d-1}\\
&&+Ck^{-d-1}k^{\frac{3}{2}d}k^{-d} \\
&\leq &\frac{C}{(k-d)^{2}}\int_{1/(k-d)}^1 j^{-d-1}(1-j)^{-d}\int_{1}^{+\infty}l^{-d-1}(l-1)^{2d-1} \dd l \dd j\\
&&+Ck^{-\frac{1}{2}d-1}\\
&\leq &C'(k-d)^{-2+d}+Ck^{-\frac{1}{2}d-1}
\end{eqnarray*}
and so the positive term has a smaller asymptotic order than the forecast error made by truncating. Therefore we have proved that in the particular case of F($d$) processes, the two prediction errors are equivalent to $Ck^{-1}$ with $C$ constant. 
\end{proof}

\begin{figure}[h]
\begin{center}
\caption{Ratio $r(k)$, $d \in ]0,1/2[$ defined in \eqref{differreur}}
\label{graphefigerreur}
\includegraphics[width=10cm]{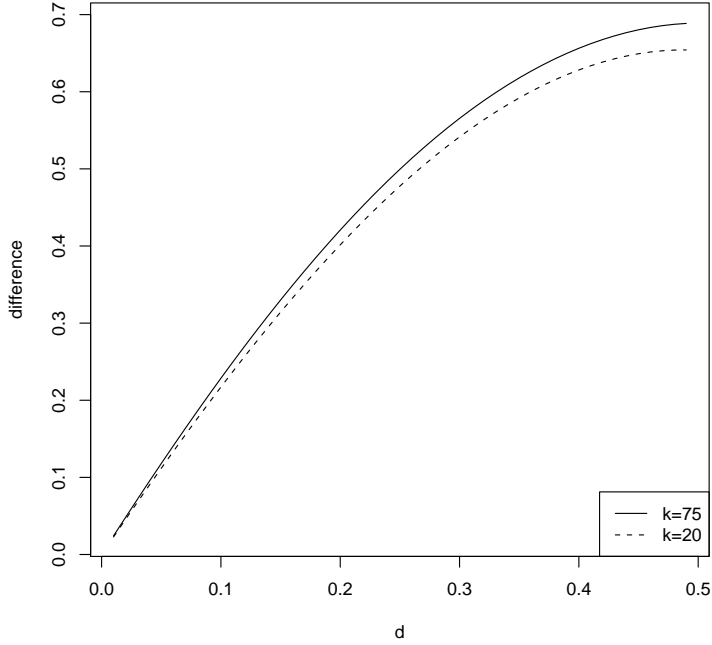}
\end{center}
\end{figure}

The two approaches to next-step prediction, by truncation to $k$ terms or by fitting an autoregressive model AR($k$) have consequently a prediction error with the same rate of convergence $k^{-1}$. So it is interesting to study how the second approach improves the prediction  The following quotient:
\begin{equation}
r(k):=\frac{\sum_{j=1}^{k}(a_{j,k}-a_j)\sum_{l=1}^{k}(a_l-a_{l,k})\sigma(j-l) 
+2\sum_{j=1}^{k}(a_{j,k}-a_j) \sum_{l=k+1}^{+\infty}a_l\sigma(j-l)}{\sum_{j=0}^{k}a_j\sum_{l=k+1}^{+\infty}a_l\sigma(j-l)} \label{differreur}
\end{equation}
is the ratio of the difference between the two prediction errors and the prediction error by truncatingn in the particular case of a fractionally integrated noise F($d$). The figure \ref{graphefigerreur} shows that the prediction by truncation incurs a larger performance loss when $d \rightarrow 1/2$. The improvement reaches 50 per cent when $d>0.3$ and $k>20$.

\subsection{Error due to Estimation of the Forecast Coefficients}

Let $(X_j)_{j \in \mathbb{Z}}$ be a stochastic process, which verifies the assumptions of section \ref{intro}, and let $(Y_j)_{j \in \mathbb{Z}}$ be a process which is independent of $(X_j)_{j \in \mathbb{Z}}$, but which has the same stochastic structure. We want to predict $X_{k+1}$ knowing $(X_j)_{j \in \llbracket 1 ,k \rrbracket}$ and we assume that forecast coefficients are estimated based on a realisation  $(Y_j)_{j \in \llbracket 1 ,T \rrbracket}$.

\par We estimate the forecast coefficients using the Yule-Walker equations \eqref{yw} where we replace the true covariances by the empirical covariances computed from the realisation $(Y_j)_{j \in \llbracket 1 ,T \rrbracket}$:
\begin{equation}
\widehat{\sigma(k)}=\frac{1}{T} \sum_{t=1}^{T-k} Y_t Y_{t+k} \label{covemp}
\end{equation}
 There exists a recursive scheme for computing the forecast coefficients. It is known as the Durbin-Levinson or innovation algorithm and it is described for example in \cite{bd}.

Let $(Y_j)_{j\in \mathbb{Z}}$ be a zero-mean process with autocovariance function $\sigma$. The coefficients $(a_{i,k})_{i \in \llbracket 1,k \rrbracket}$ satisfy the Yule-Walker $k^{\mathrm{th}}$ equations:
\begin{displaymath}
\forall j \in \llbracket 1,k \rrbracket, \quad \sigma(j)=\sum_{u=1}^k \sigma(u-j)a_{u,k}.
\end{displaymath}
If we let $v(0)=\sigma(0)$ and $a_{1,1}=\sigma(1)/\sigma(0)$, then we have for any integer $n$:
\begin{eqnarray*}
a_{n,n}&=&\big[\sigma(n)-\sum_{j=1}^{n-1}a_{j,n-1}\sigma(n-j) \big]\frac{1}{v(n-1)}\\
\left(\begin{array}{c}
a_{1,n}\\
\vdots\\
a_{n-1,n}
\end{array}\right)
&=&
\left(\begin{array}{c}
a_{1,n-1}\\
\vdots\\
a_{n-1,n-1}
\end{array}\right)
-a_{n,n}
\left(\begin{array}{c}
a_{n-1,n-1}\\
\vdots\\
a_{1,n-1}
\end{array}\right)\\
v(n)&=&v(n-1)(1-a_{n,n}^2).
\end{eqnarray*}

We denote by $(\widehat{a_{1,k}},\ldots,\widehat{a_{k,k}})$ the respective solutions to the Yule-Walker equations obtained by replacing the covariances by theirs estimates the empirical covariances defined in \eqref{covemp}. Contrary to Section \ref{parametrique}, the estimation of the forecast coefficients is non-parametric.
\par Another way to estimate the coefficients has been considered by \cite{dist}. Our method borrow the idea (see section \ref{yaj}) that the coefficients of the AR($k$) minimize: 
\begin{equation}
\int_{-\pi}^{\pi}f(\lambda)/g(\lambda,\Phi)\dd \lambda \label{distmodele}.
\end{equation}
If we replace in \eqref{distmodele} the spectral density by the periodogram $I_T$ :
\begin{displaymath}
I_T(\lambda)=\frac{\left| \sum_{t=1}^T X_t e^{i t \lambda} \right|^2}{2 \pi T},
\end{displaymath}
then:
\begin{equation}
(\widehat{a_{1,k}},\ldots,\widehat{a_{k,k}})= \underset{\Phi}{\mathrm{argmin}}\int_{-\pi}^{\pi}I_T(\lambda)/g(\lambda,\Phi)\dd \lambda\label{min}.
\end{equation}

\par From now on we incorporate the effects of estimation of the AR($k$) coefficients using a realisation of length $T$, as $T \rightarrow + \infty$ and study the mean-squared prediction error due to this estimation.
We define $\widehat{X_{T,k}}(1)$ the predictor with all the coefficients $a_{j,k}$ replaced by their estimates:
\begin{displaymath}
\widehat{X_{T,k}}(1):= \sum_{j=1}^k \widehat{a_{j,k}} X_{k+1-j}
\end{displaymath}

More precisely, we study the mean-squared difference between the predictor with the estimated coefficients $\widehat{a_{j,k}}$ and the predictor with the true coefficients $a_{j,k}$:\\
\begin{eqnarray*}
&&\E\left[\left(\widehat{X_{T,k}}(1)- \widehat{X_{k}}(1)\right)^2  \right] \\
&=&\mathbb{E}\left( \left( \widehat{a_{1,k}}-a_{1,k},\ldots,\widehat{a_{k,k}}-a_{k,k}\right) 
\left( \begin{array}{c}
X_{k} \\
\vdots \\
X_1
\end{array}\right) 
\right) ^2 \nonumber \\
&=& \mathrm{trace}\left(\mathbb{E} \left(
\left( 
\begin{array}{c}
\widehat{a_{1,k}}-a_{1,k} \\
\vdots \\
\widehat{a_{k,k}}-a_{k,k}
\end{array}
\right) 
\left( \widehat{a_{1,k}}-a_{1,k},\ldots,\widehat{a_{k,k}}-a_{k,k}\right)
\right)
\mathbb{E} \left(
\left( \begin{array}{c}
X_{k} \\
\vdots \\
X_1
\end{array}\right)( X_{k},\ldots ,X_1)\right) \right) \\
&=&\mathrm{trace}\left(\mathbb{E} \left(
\left( 
\begin{array}{c}
\widehat{a_{1,k}}-a_{1,k} \\
\vdots \\
\widehat{a_{k,k}}-a_{k,k}
\end{array}
\right) \left(\widehat{a_{1,k}}-a_{1,k},\ldots,\widehat{a_{k,k}}-a_{k,k} \right)
\right)
\Sigma_k \right).
\end{eqnarray*}
First we estimate the covariance matrix:
\begin{displaymath}
\mathbb{E} \left(
\left( 
\begin{array}{c}
\widehat{a_{1,k}}-a_{1,k} \\
\vdots \\
\widehat{a_{k,k}}-a_{k,k}
\end{array}
\right) \left(\widehat{a_{1,k}}-a_{1,k},\ldots,\widehat{a_{k,k}}-a_{k,k} \right)
\right).
\end{displaymath}
For later convenience, we now introduce the vector:
\begin{displaymath}
\textbf{1}_k^*:=\underbrace{(1, \ldots,1)}_{k}
\end{displaymath}
and the $(k \times k)$ matrix:
\begin{displaymath}
\textbf{1}_{k,k}:=\textbf{1}_k \textbf{1}_k^*.
\end{displaymath}
We now state the theorem which allows us to conclude.

\begin{thm} \label{passagecovcoeff}
We assume that the process $(Y_n)_{n \in \mathbb{Z}}$ is Gaussian, that its autocovariance function $\sigma$ verifies:
 \begin{displaymath}
 \sigma(j) \sim \lambda j^{2d-1} \textrm{ with } \lambda>0,
 \end{displaymath}
that the coefficients of its infinite moving average representation $b_j$ verify:
$$b_j \sim  \delta j^{d-1}\textrm{ with } \delta>0,$$
and finally that the white noise process $(\e_n)_{n \in \mathbb{Z}}$ is such that $\forall n \in \mathbb{Z},\: \E(\e_n^4)<+\infty$.
We will denote by $g_{i,j}$ the function:
\begin{displaymath}
\begin{array}{cccc}
g_{i,j}: & \mathbb{R}^{k+1}& \rightarrow & \mathbb{R}\\
 & (x_0,\ldots,x_k)& \mapsto &(y_i-a_{i,k})(y_j-a_{j,k})
\end{array}
\end{displaymath}
where
\begin{displaymath}
\left( \begin{array}{c}
y_1 \\
\vdots \\
y_k
\end{array}\right)=
\left( \begin{array}{cccc}
x_0&x_1& \ldots & x_{k} \\
x_1&x_0& \ddots & x_{k-1} \\
\vdots& \ddots& \ddots \vdots \\
x_{k}&x_{k-1}&\ldots &x_0
\end{array}\right) ^{-1}
\left( \begin{array}{c}
x_1 \\
\vdots \\
x_k
\end{array}\right) .
\end{displaymath}
Then
\begin{eqnarray*}
&&\mathbb{E}(g_{i,j}(\widehat{\sigma(0)},\widehat{\sigma(1)},\ldots,\widehat{\sigma(k)}))\\&&=
\begin{cases}
\left(1-\sum_{r=1}^k a_{r,k} \right)^2C n^{4d-2}\left(\Sigma_k^{-1}\textbf{1}_{k,k} \Sigma_k^{-1}\right)_{(i,j)}+\mathrm{O}(n^{6d-3})\quad \textrm{if }\frac{1}{4}<d<\frac{1}{2}\\
\left(1-\sum_{r=1}^k a_{r,k} \right)^2D \frac{\ln(n)}{n}\left(\Sigma_k^{-1}\textbf{1}_{k,k} \Sigma_k^{-1}\right)_{(i,j)}+\mathrm{O}(n^{-3/2}) \quad \textrm{if }d=1/4 \\
n^{-1} 4\left( \Sigma_k^{-1} H \Sigma_k^{-1} \right)_{(i,j)}+\mathrm{O}(n^{-3/2}) \quad \textrm{if }0<d<\frac{1}{4}
\end{cases}
\end{eqnarray*}
where $C$ and $D$ are constants independent of $n$ and $k$. The definition of the matrix $H$ follows. We define $h$ as $h(\lambda) = |1- \sum_{r=1}^k a_{r,k} \mathrm{e}^{ir\lambda} |^2$ and we denote by $h^{(r)}$ the derivative of the function $h$ with respect to $ a_{r,k}$. The $(i,j)$-th entry of the matrix $H$ is given by:
\begin{displaymath}
H_{i,j}:=\int_{-\pi}^{\pi} h^{(i)}(\lambda) h^{(j)}(\lambda) f^2(\lambda) \dd \lambda.
\end{displaymath} \label{thmdl}
\end{thm}
\begin{proof}
The proof is given in Appendix.
\end{proof}

Next we estimate the asymptotic behaviour of $\E\left[\left(\widehat{X_{T,k}}(1)- \widehat{X_{k}}(1)\right)^2  \right] $ and we state the following Theorem which gives an estimation of the mean-squared error when $d\geq 1/4$:
\begin{thm} We assume that the assumptions of the Theorem \ref{passagecovcoeff} hold. \label{d>1/4}
We assume also that the spectral density of the process is such that:
\begin{displaymath}
\forall x \in [-\pi, \pi],\quad f(x)=f_{d}(x)L(x)
\end{displaymath}
with $f_{d}$ defined by:
\begin{displaymath}
\forall x \in [-\pi, \pi],\quad f_{d}(x)=2^{-2d-1}\pi^{-1}\left( \sin^2(x/2)\right)^{-d} 
\end{displaymath} 
and $L$ a positive, integrable on $[-\pi, \pi]$, continuous at 0 and bounded below by a positive constant. \\
If $d = 1/4$ then
\begin{displaymath}
\E\left[\left(\widehat{X_{T,k}}(1)- \widehat{X_{k}}(1)\right)^2  \right]  =\OO\left( \frac{\log(T)\sqrt{k}}{T}\right)  
\end{displaymath}  \\
and if $d \in \left]1/4,1/2\right[ $, we thus get
\begin{displaymath}
\E\left[\left(\widehat{X_{T,k}}(1)- \widehat{X_{k}}(1)\right)^2  \right] =\OO\left(\frac{k^{1-2d}}{T^{2-4d}}\right).
\end{displaymath}
\end{thm}

\paragraph*{Remark} The assumption that $L$ and so $f$ are bounded below by a positive constant is not a new very restrictive assumption. Since we have assumed that the process admits an infinite autoregressive representation:
\begin{displaymath}
\e_n=\sum_{j=0}^{\infty} a_j X_{n-j}, 
\end{displaymath}
where the coefficients $a_j$ are absolutely summable, we have that the spectral density can be written as:
\begin{displaymath}
f(\lambda)=\frac{1}{\left\vert \sum_{j=0}^{\infty}a_j \mathrm{e}^{ij\lambda}\right \vert^2}
\end{displaymath}
 and consequently the spectral density can not vanish on $[-\pi,\pi[$.

\begin{proof}
Applying the last theorem, we obtain that if $d = 1/4$ then 
\begin{displaymath}
\E\left[\left(\widehat{X_{T,k}}(1)- \widehat{X_{k}}(1)\right)^2  \right] =\OO\left( \textrm{trace}\left( \frac{\log(T)}{T} \left(\sum_{j=0}^{k} a_{j,k}\right)^2   \Sigma_k^{-1}
\textbf{1}_{k,k} 
\right)\right)
\end{displaymath} 
and if $d \in \left]1/4,1/2\right[ $ then
\begin{displaymath}
\E\left[\left(\widehat{X_{T,k}}(1)- \widehat{X_{k}}(1)\right)^2  \right] =\OO\left(  \textrm{trace}\left( \frac{1}{T^{2-4d}} \left(\sum_{j=0}^{k} a_{j,k}\right)^2   \Sigma_k^{-1}\textbf{1}_{k,k} 
\right)\right).
\end{displaymath}

 First we estimate $\left| \sum_{j=0}^{k} a_{j,k}\right|$. We write this like:
 \begin{eqnarray}
\left| \sum_{j=0}^{k} a_{j,k}\right| &\leq&  \sqrt{k} \sqrt{\sum_{j=1}^{k}|a_{j,k}-a_j|}+\sum_{j=0}^{k}|a_j|\nonumber.
\end{eqnarray}
We follow the proof of Theorem 3.3 of \cite{inoue} about the convergence of the sequence of the misspecified AR($k$) model coefficients to the R($\infty$) representation coefficients. We shall remark that there exists $C_1$, $C_2$ and $K$ such that if:
\begin{displaymath}
k\geq K, \quad k\left( a_{j,k}-a_j\right)\leq C_1 \sum_{u=k-j}^{+\infty} |a_u| + C_2 \sum_{u=j}^{+\infty} |a_u|
\end{displaymath}
and $C$ is a generic constant:
\begin{displaymath}
\textrm{si}\quad k\geq K, \quad k\left( a_{j,k}-a_j\right)\leq C\left(  \sum_{u=k-j}^{+\infty} |a_u| + \sum_{u=j}^{+\infty} |a_u|\right) 
\end{displaymath}
We thus get that if $k\geq K$:
\begin{eqnarray}
\sum_{j=1}^{k}(a_{j,k}-a_j)^2&\leq&\frac{C}{k^2}\sum_{j=1}^k  \left[ \left( \sum_{u=k-j}^{+\infty} |a_u| \right) ^2+\left( \sum_{u=j}^{+\infty} |a_u| \right)\right]  ^2\nonumber \\
&\leq&  \frac{C}{k^2}\sum_{j=1}^k (k-j+1)^{-2d}+j^{-2d}\nonumber \\
&\leq&  \OO\left( k^{-2d-1}\right) .
\end{eqnarray}
So we may conclude that:
\begin{eqnarray}
\left| \sum_{j=0}^{k} a_{j,k}\right| &\leq& \OO(k^{-d})+\sum_{j=0}^{+\infty}|a_j|\nonumber\\
&=& \OO(1). \label{maj}
\end{eqnarray}
 Next we have to study the asymptotic properties of:
 \begin{displaymath}
 \mathrm{trace} \left( \Sigma_k^{-1}
\textbf{1}_{k,k} \right) =\begin{array}{ccc}
\left( 1 \ldots 1 \right) 
\end{array}
\Sigma_k^{-1}
\left( 
\begin{array}{c}
1 \\ \vdots \\ 1 
\end{array} \right)
 \end{displaymath}

Then applying Theorem 6.1 of \cite{adenstedt} under the assumptions of theorem \ref{d>1/4} we obtain the following asymptotic equivalent:
\begin{displaymath}
\begin{array}{ccc}
\left( 1 \ldots 1 \right) 
\end{array}
\Sigma_k^{-1}
\left( 
\begin{array}{c}
1 \\ \vdots \\ 1 
\end{array} \right) \sim\left( \frac{k^{1-2d}  \Gamma(-2d+1) L(0)}{\beta(-d+1,-d+1)} \right)^{-1}
\end{displaymath}
where $\Gamma$ and $\beta$ are respectively the gamma function and the beta function. The result follows.

\end{proof}

The last case is when  $0<d<1/4$:
\begin{thm}We assume that the assumptions of \ref{passagecovcoeff} hold. We assume also that the spectral density $f$ is bounded above by a constant positive.
If $0<d<1/4$ then
\begin{displaymath}
\E\left[\left(\widehat{X_{T,k}}(1)- \widehat{X_{k}}(1)\right)^2  \right] =\OO\left( \frac{k}{T}\right) .
\end{displaymath}
\end{thm}

\begin{proof}
We call $(\Phi_i:=\Phi_{i,1}+\ldots+\Phi_{i,i}x^{i-1})_{i \in \mathbb{N^*}}$ the orthonormal polynomials associated with the spectral density $f$, that is to say $\Phi_i$ is a $(i-1)^{\mathrm{th}}$ degree polynomial such that
\begin{displaymath}
\forall j , l \in \mathbb{N^*}, \quad \int_{-\pi}^{\pi} f(\lambda)\Phi_j\left(\mathrm{e}^{i\lambda}\right) \Phi_l\left(\mathrm{e}^{-i\lambda}\right) \dd \lambda = \delta_{j,l}
\end{displaymath}
where $\delta$ is the Kronecker delta.
We then define the matrix $T_k$ by:
\begin{displaymath}
T_k=
\left( 
\begin{array}{cccc}
\Phi_{1,1} &0 & \ldots & 0\\
\Phi_{2,1}& \Phi_{2,2}& \ddots& 0\\
\vdots & \vdots & \ddots & \vdots \\
\Phi_{k,1} & \Phi_{k,2}&\ldots &\Phi_{k,k}
\end{array}
\right) 
\end{displaymath}
$T_k$ verifies the following conditions:
\begin{displaymath}
T_k \Sigma_k T_k^* = Id_k 
\end{displaymath}
and so
\begin{equation}
\Sigma_k^{-1}=T_k^*T_k. \label{Tk}
\end{equation}

Using \eqref{Tk}, we obtain that:
\begin{eqnarray}
\E\left[\left(\widehat{X_{T,k}}(1)- \widehat{X_{k}}(1)\right)^2  \right] &=&\frac{1}{T} \mathrm{trace} \left( \Sigma_k^{-1} H\right) \nonumber\\
 &=&\frac{1}{T}\mathrm{trace} \left(T_k H T_k^*\right) \label{traceamaitriser}
 \end{eqnarray}
with $H$ defined in Theorem \ref{thmdl}. 
We define $G_k:\lambda \mapsto \sum_{j=0}^k a_{j,k} \mathrm{e}^{ij \lambda}$.\\$ \E\left[\left(\widehat{X_{T,k}}(1)- \widehat{X_{k}}(1)\right)^2  \right] $ can therefore be rewritten like:
\begin{eqnarray*}
&&\frac{1}{T}\mathrm{trace}\left(\left( \int_{-\pi}^{\pi} f^2(\lambda) \R\left(  G_k(\lambda)\Phi_j(\mathrm{e}^{ij\lambda})\right) \R\left(  G_k(-\lambda)\Phi_l(\mathrm{e}^{-il\lambda})\right) \dd \lambda \right)_{j,l \in \llbracket 1,k \rrbracket}\right) \\
&=&\frac{1}{T} \mathrm{trace}\left(\left( \int_{-\pi}^{\pi} f^2(\lambda) \R\left(  G_k(\lambda)\Phi_j(\mathrm{e}^{ij\lambda})  G_k(-\lambda)\Phi_l(\mathrm{e}^{-il\lambda})\right) \dd \lambda \right)_{j,l \in \llbracket 1,k \rrbracket}\right) \\
&+&\frac{1}{T}\mathrm{trace}\left(\left( \int_{-\pi}^{\pi} f^2(\lambda) \textrm{Im}\left(  G_k(\lambda)\Phi_j(\mathrm{e}^{ij\lambda})\right) \textrm{Im}\left(  G_k(-\lambda)\Phi_l(\mathrm{e}^{-il\lambda})\right) \dd \lambda \right)_{j,l \in \llbracket 1,k \rrbracket}\right) 
\end{eqnarray*}
because $\R(ab)=\R(a)\R(b)-\textrm{Im}(a)\textrm{Im}(b)$.
For later convenience, we note:
\begin{eqnarray*}
A&:=&\left( \int_{-\pi}^{\pi} f^2(\lambda) \R\left(  G_k(\lambda)\Phi_j(\mathrm{e}^{i\lambda})\right) \R\left(  G_k(-\lambda)\Phi_l(\mathrm{e}^{-i\lambda})\right) \dd \lambda \right)_{j,l \in \llbracket 1,k \rrbracket} \\
B&:=&\left( \int_{-\pi}^{\pi} f^2(\lambda) \R\left(  G_k(\lambda)\Phi_j(\mathrm{e}^{i\lambda})  G_k(-\lambda)\Phi_l(\mathrm{e}^{-i\lambda})\right) \dd \lambda \right)_{j,l \in \llbracket 1,k \rrbracket} \\
&=&\left( \int_{-\pi}^{\pi} f^2(\lambda)  \vert G_k(\lambda)\vert^2 \Phi_j(\mathrm{e}^{i\lambda})  \Phi_l(\mathrm{e}^{-i\lambda}) \dd \lambda \right)_{j,l \in \llbracket 1,k \rrbracket} \\
C&:=&\left( \int_{-\pi}^{\pi} f^2(\lambda) \textrm{Im}\left(  G_k(\lambda)\Phi_j(\mathrm{e}^{i\lambda})\right) \textrm{Im}\left(  G_k(-\lambda)\Phi_l(\mathrm{e}^{-i\lambda})\right) \dd \lambda \right)_{j,l \in \llbracket 1,k \rrbracket}
\end{eqnarray*}
Then we have $A=B+C$. We will prove that $A$, $B$ and $-C$ are symmetric and positive matrices, which implies that $0 \leq \mathrm{trace}(A) \leq \mathrm{trace}(B)$. First we study the symmetry: $A$ is symmetric because the real part of a complex is equal to that of its conjugate, $B$ is symmetric because $\lambda \mapsto f^2(\lambda)\vert G_k(\lambda) \vert^2$ is a symmetric function and $C$ is symmetric because the imaginary part is equal to the negative of the imaginary part of its conjugate. Next we study the positivity. Let $q:=(q_1, \ldots ,q_k)$ be a vector. We have:
\begin{eqnarray*}
q A q^*&=& \int_{-\pi}^{\pi} f^2(\lambda) \R\left( \sum_{j=1}^k G_k(\lambda) q_j\Phi_j(\mathrm{e}^{i\lambda})\right) \R\left( \sum_{l=1}^k  G_k(-\lambda)q_l\Phi_l(\mathrm{e}^{-i\lambda})\right) \dd \lambda 
\geq 0 \\
q B q^* &=& \int_{-\pi}^{\pi} f^2(\lambda)  \vert G_k(\lambda)\vert^2\sum_{j=1}^k  q_j\Phi_j(\mathrm{e}^{i\lambda}) \sum_{l=1}^k  q_l\Phi_l(\mathrm{e}^{-i\lambda}) \dd \lambda 
\geq 0 \\
q Cq^* &=& \int_{-\pi}^{\pi} f^2(\lambda) \textrm{Im}\left( \sum_{j=1}^k  G_k(\lambda)q_j\Phi_j(\mathrm{e}^{i\lambda})\right) \textrm{Im}\left( \sum_{j=1}^k  G_k(-\lambda)q_l\Phi_l(\mathrm{e}^{-i\lambda})\right) \dd \lambda\leq0.
\end{eqnarray*}

The traces of these matrices $A$ , $B$ et $-C$ are equal to the sum of theirs eigenvalues since they are symmetric and thus diagonalizable. Because these matrices are positive, all their eigenvalues are positive and the traces are also positive. Therefore we obtain that:
\begin{equation}
0 \leq \mathrm{trace}(A) \leq \mathrm{trace}(B). \label{trace}
\end{equation}
To find a bound for $\mathrm{trace}(A)$, it is sufficient to find a bound for $\mathrm{trace}(B)$:
\begin{eqnarray*}
\mathrm{trace}(B) &=& \sum_{j=1}^{k} \int_{-\pi}^{\pi} f^2(\lambda)  \vert G_k(\lambda)\vert^2 \Phi_j(\mathrm{e}^{i\lambda})  \Phi_j(\mathrm{e}^{-i\lambda}) \dd \lambda\\
&=& \int_{-\pi}^{\pi} f^2(\lambda)  \vert G_k(\lambda)\vert^2  K_k(\mathrm{e}^{i\lambda},\mathrm{e}^{i\lambda})\dd \lambda
\end{eqnarray*}
where $K_k$ is the reproducing kernel defined by:
\begin{displaymath}
\forall x,y \in \mathbb{C}, K_k(x,y) =\sum_{j=1}^{k} \Phi_j(x) \Phi_j(\bar{y}).
\end{displaymath}
We have assumed that the spectral density $f$ is bounded from below by a positive constant $c$, so we can apply the Theorem 2.2.4 of \cite{simon} and we get:
\begin{displaymath}
\forall \lambda \in [-\pi, \pi], K_k(\mathrm{e}^{i\lambda},\mathrm{e}^{i\lambda}) \leq k \frac{2\pi}{c}.
\end{displaymath}
We look for a bound for $\vert G_k(\lambda)\vert^2 $:
\begin{eqnarray*}
\forall \lambda \in [-\pi, \pi],\vert G_k(\lambda)\vert^2 &\leq& \left(\sum_{j=0}^k |a_{j,k}| \right) ^2\\
&=&\OO (1)
\end{eqnarray*}
as we have proven in \eqref{maj}. This bound is independent of $\lambda$. We finally notice that if $0<d< \frac{1}{4}$ then $f$ is square integrable.
So we obtain that:
\begin{displaymath}
\tr(B) = \OO(k) 
\end{displaymath}
and we conclude using \eqref{traceamaitriser} and \eqref{trace} that:
\begin{displaymath}
\E\left[\left(\widehat{X_{T,k}}(1)- \widehat{X_{k}}(1)\right)^2  \right] = \OO \left(\frac{k}{T} \right). 
\end{displaymath}


\end{proof}

\subsection{Conclusion}

\par Fitting an AR($k$) model also involves two mean-squared error components: the first is due to fitting a misspecified model and is bounded by $\OO(k^{-1})$ and the second is due to the estimation of the Yule-Walker coefficients $a_{j,k}$ from a independent realisation of length $T$ and is bounded by $\OO(k/T)$ if $0<d<1/4$ (\cite{bhansali} has the same asymptotic equivalent for short memory processes), bounded by $\OO(k^{1/2}\log(T)/T)$ if $d=1/4$ and bounded by $\OO(k^{1-2d}/T^{2-4d})$ if $1/4<d<1/2$.
As in Section \ref{conclusion1}, if we want to compare the two types of forecast error, we need to state a relation between $k$ and $T$ and moreover distinguish 3 cases for the value of $d$.

\paragraph*{} In both methods by truncating to $k$ terms the Wiener-Kolmogorov predictor or by fitting an AR($k$) model, the mean-squared error of prediction due to the method is bounded by $\OO(k^{-1})$. Nevertheless, the factor of $k^{-1}$ in this equivalent depends on $d$. We have shown that the factor tends to infinity when $d$ tends to $1/2$ in the method by truncation in the special case of fractionally integrated noise (Section \ref{erreurtronc}) so that the error increases for $d$ near  $1/2$. Moreover for this value of $d$, figure \eqref{graphefigerreur} show that fitting an AR($k$) model greatly reduces the error.
For the errors due to the estimation of the forecast coefficients, the method by truncation is optimal since if we assume that $T/k$ tends to infinity (necessary condition to have some mean-squared error which converges to 0), then for all $d$ in $]0,1/2[\backslash\{1/4\}$:
\begin{displaymath}
\E\left( \widetilde{X'_{T,k}}(1)-\widetilde{X'_k}(1)
\right)^2=\mathrm{o}\left(\E\left[\left(\widehat{X_{T,k}}(1)- \widehat{X_{k}}(1)\right)^2  \right] \right) .
\end{displaymath}
In the end, we have so to consider the value of long-memory parameter $d$, the length of the series $k$ and $T$ to decide on a prediction method.
\section{Appendix}

\subsection{Proof of Lemma \ref{sommedoublecomp} } \label{annexe1}
Let $g$ be the function $(l,j) \mapsto j^{-d-1+\delta}l^{-d-1+\delta}\vert l-j\vert^{2d-1+\delta}$. Let $m$ and $n$ be two integers. We assume that $\delta<1-2d $ and that $m \geq \frac{\delta -d-1}{\delta+2d-1}$ for all $\delta \in \left] 0,\frac{\delta -d-1}{\delta+2d-1}\right[ $. 
We introduce $A_{n,m}$ the square $[n,n+1]\times[m,m+1]$. If  $n \geq m+1$ then
\begin{displaymath}
\int_{A_{n,m}}g(l,j) \dd j \dd l \geq g(n+1,m).
\end{displaymath}
\begin{proof}
We restrict the domain of $g$ to the square $A_{n,m}$. First we will show that $g(.,j)$ is a decreasing and then we compute its derivative:
\begin{eqnarray*}
\left( g(j,.)\right)'(l)&=&\left[ (-d-1+\delta) l^{-1}+(2d-1+\delta)(l-j)^{-1} \right] j^{-d-1+\delta}l^{-d-1+\delta}( l-j)^{2d-1+\delta}\\
&\leq& 0
\end{eqnarray*}
since $\delta<1-2d $.
We show then that $g(l,.)$ is increasing:
\begin{eqnarray*}
\left( g(.,l)\right)'(j)&=&\left[(-d-1+\delta)j^{-1}-(2d-1+\delta)(l-j)^{-1}\right] j^{-d-1+\delta}l^{-d-1+\delta}( l-j)^{2d-1+\delta}\\
& \geq&0
\end{eqnarray*}
because
\begin{displaymath}
j\geq \frac{\delta -d-1}{\delta+2d-1}.
\end{displaymath}
Then the function $g$ attains its minimum at $(n+1,m)$ and we have 
\begin{eqnarray*}
\forall(l,j) \in A_{n,m},g(l,j) &\geq &g(n+1,m)\\
 \int_{A_{n,m}}g(l,j) \dd j \dd l &\geq &\int_{A_{n,m}}g(n+1,m)\dd j \dd l\\
 \int_{A_{n,m}}g(l,j) \dd j \dd l &\geq &g(n+1,m).
\end{eqnarray*}
The results follows.
 
\end{proof}

\subsection{Proof of Theorem \ref{passagecovcoeff}}
By assumption, the process $(Y_n)_{n \in \mathbb{Z}}$ is Gaussian. We also assume that its autocovariance function $\sigma$ verifies:
 \begin{displaymath}
 \sigma(j) \sim \lambda j^{2d-1} \textrm{ with } \lambda>0,
 \end{displaymath}
that the coefficients of its moving average representation $b_j$ are such that:
$$b_j \sim  \delta j^{d-1}\textrm{ with } \delta>0,$$
and that the white-noise series $(\e_n)_{n \in \mathbb{Z}}$ has finite fourth moments.
 Let $g_{i,j}$ be the function:
\begin{displaymath}
\begin{array}{cccc}
g_{i,j}: & \mathbb{R}^{k+1}& \rightarrow & \mathbb{R}\\
 & (x_0,\ldots,x_k)& \mapsto &(y_i-a_{i,k})(y_j-a_{j,k})
\end{array}
\end{displaymath}
with
\begin{displaymath}
\left( \begin{array}{c}
y_1 \\
\vdots \\
y_k
\end{array}\right)=
\left( \begin{array}{cccc}
x_0&x_1& \ldots & x_{k} \\
x_1&x_0& \ddots & x_{k-1} \\
\vdots& \ddots& \ddots &\vdots \\
x_{k}&x_{k-1}&\ldots &x_0
\end{array}\right) ^{-1}
\left( \begin{array}{c}
x_1 \\
\vdots \\
x_k
\end{array}\right) .
\end{displaymath}
Therefore
\begin{eqnarray*}
&&\mathbb{E}(g_{i,j}(\widehat{\sigma(0)},\widehat{\sigma(1)},\ldots,\widehat{\sigma(k)}))\\&&=
\begin{cases}
\left(1-\sum_{r=1}^k a_{r,k} \right)^2C n^{4d-2}\left(\Sigma_k^{-1}\textbf{1}_{k,k} \Sigma_k^{-1}\right)_{(i,j)}+\mathrm{O}(n^{6d-3})\quad \textrm{if }\frac{1}{4}<d<\frac{1}{2}\\
\left(1-\sum_{r=1}^k a_{r,k} \right)^2D \frac{\ln(n)}{n}\left(\Sigma_k^{-1}\textbf{1}_{k,k} \Sigma_k^{-1}\right)_{(i,j)}+\mathrm{O}(n^{-3/2}) \quad \textrm{if }d=1/4 \\
n^{-1} 4\left( \Sigma_k^{-1} H \Sigma_k^{-1} \right)_{(i,j)}+\mathrm{O}(n^{-3/2}) \quad \textrm{if }0<d<\frac{1}{4}
\end{cases}
\end{eqnarray*}
where $C$ and $D$ are constants and independent of $n$ and $k$. The definition of matrix $H$ follows. We define $h$ as $h(\lambda) = |1- \sum_{r=1}^k a_{r,k} \mathrm{e}^{ir\lambda} |^2$ and we denote by $h^{(r)}$ the derivative of the function $h$ with respect to $ a_{r,k}$. The $(i,j)$-th entry of the matrix $H$ is given by:
\begin{equation}
H_{i,j}:=\int_{-\pi}^{\pi} h^{(i)}(\lambda) h^{(j)}(\lambda) f^2(\lambda) \dd \lambda \label{H}
\end{equation}
\begin{proof}
We write a $2^{\mathrm{nd}}$ order Taylor expansion of the function $g_{i,j}$ applying Theorem 5.4.3 in \cite{dl} as in the Section \ref{parametrique}. We will refer to the following version.:
 \begin{enumerate}[(i)]
 \item If $\E\left( |\widehat{\sigma(k)}-\sigma(k)|^3\right) = \OO(a_n)$;
 \item if $g_{i,j}$ is uniformly bounded;
 \item if the first and the second derivatives of $g_{i,j}$ are continuous and bounded functions on a neighbourhood of $(\sigma(0),\ldots,\sigma(k))$
 \end{enumerate}
then
 \begin{eqnarray*}
 &&\E\left( g_{i,j}(\widehat{\sigma(0)},\widehat{\sigma(1)},\ldots,\widehat{\sigma(k)})\right)\\ &=&g_{i,j}(\sigma(0),\ldots,\sigma(k))+\sum_{l=0}^{k} \E(\widehat{\sigma(l)}-\sigma(l))\frac{\partial g_{i,j}}{\partial x_l}(\sigma(0),\ldots,\sigma(k))\\
 &&+\frac{1}{2}\sum_{l=0}^{k}\sum_{m=0}^{k}\frac{\partial^2 g_{i,j} }{\partial x_l \partial x_m}(\sigma(0),\ldots,\sigma(k))\E\left( (\widehat{\sigma(l)}-\sigma(l))(\widehat{\sigma(m)}-\sigma(m))\right) +\OO(a_n).
 \end{eqnarray*}
 We first verify that the assumptions hold.
 We need a bound for the third order moments of the empirical covariances.
\begin{lm}\label{lmmoments}
\begin{equation}
\mathbb{E}\left[ \Big| \frac{1}{n}\sum_{t=0}^{n-k}X_tX_{t+k}-\sigma(k) \Big|^3 \right] =
\begin{cases}
\OO(n^{-3/2})& \textrm{ if $d \leq 1/4$} \\
\OO(n^{6d-3})&\textrm{ if $d > 1/4$}
\end{cases}
\end{equation}
\end{lm}
\begin{proof}
Lemma \ref{lmmoments} is proven in Section \ref{lemme2}.
\end{proof}
In this way we obtain a bound for the rest of the Taylor series.
Moreover $g_{i,j}$ is an uniformly bounded function since its results are the coefficients of the autoregressive process. Since the derivatives of $g_{i,j}$ at $(\sigma(0),\ldots,\sigma(k))$ are finite, there exits a neighbourhood of $(\sigma(0),\ldots,\sigma(k))$ such that on this all the derivatives are uniformly bounded. So we apply the Theorem 5.4.3 in \cite{dl}.
First we note that:
\begin{displaymath}
g_{i,j}(\sigma(0),\ldots,\sigma(k))=0
\end{displaymath}
and 
\begin{eqnarray}
\forall l \in \llbracket 0,k \rrbracket, \quad \frac{\partial g_{i,j}}{\partial x_l}(\sigma(0),\ldots,\sigma(k))&=&\frac{\partial y_i}{\partial x_l}(\sigma(0),\ldots,\sigma(k))(y_j-a_{j,k})\nonumber \\
&&+(y_i-a_{i,k})\frac{\partial y_j}{\partial x_l}(\sigma(0),\ldots,\sigma(k))\nonumber \\
\frac{\partial g_{i,j}}{\partial x_l}(\sigma(0),\ldots,\sigma(k))&=&0 \label{annule}
\end{eqnarray}
because
\begin{displaymath}
\forall i \in \llbracket 1,k \rrbracket, \quad  y_i(\sigma(0),\ldots,\sigma(k))-a_{i,k}=0.
\end{displaymath}
From the Taylor series and Lemma \ref{lmmoments}, it follows that:
\begin{eqnarray*}
&&\mathbb{E}(g_{i,j}(\widehat{\sigma(0)},\widehat{\sigma(1)},\ldots,\widehat{\sigma(k)}))\\&&= \begin{cases}
\sum_{l=0}^{k}\sum_{m=0}^{k}\frac{\partial g_{i,j} }{\partial x_l \partial x_m}(\sigma(0),\ldots,\sigma(k)) \E((\widehat{\sigma(l)}-\sigma(l))(\widehat{\sigma(m)}-\sigma(m)))+\OO(n^{-3/2})\\
 \textrm{if } 0<d \leq 1/4\\
 \sum_{l=0}^{k}\sum_{m=0}^{k}\frac{\partial g_{i,j} }{\partial x_l \partial x_m}(\sigma(0),\ldots,\sigma(k)) \E((\widehat{\sigma(l)}-\sigma(l))(\widehat{\sigma(m)}-\sigma(m)))+\OO(n^{6d-3}) \\
 \textrm{if } 1/4<d<1/2 \end{cases}
\end{eqnarray*}
According to the results from \cite{convcov}, we shall compute the second term of the Taylor series.
First we note that:
\begin{displaymath}
\frac{\partial g_{i,j} }{\partial x_l \partial x_m}(\sigma(0),\ldots,\sigma(k))=\frac{\partial y_i}{\partial x_l}\frac{\partial y_j}{\partial x_m}(\sigma(0),\ldots,\sigma(k))+\frac{\partial y_i}{\partial x_m}\frac{\partial y_j}{\partial x_l}(\sigma(0),\ldots,\sigma(k))
\end{displaymath}
because the other terms vanish at $(\sigma(0),\ldots,\sigma(k))$ by \eqref{annule}. Moreover we can apply \cite{convcov}:
\begin{eqnarray*}
&&\E((\widehat{\sigma(l)}-\sigma(l))(\widehat{\sigma(m)}-\sigma(m)))\\ && \stackrel{n \rightarrow + \infty}{\sim}
\begin{cases} C n^{4d-2} &  \textrm{if} \quad\frac{1}{4}<d<\frac{1}{2} \\
D n^{-1} \ln(n)& \textrm{if}\quad d=\frac{1}{4} \\
 n^{-1}\left( \sum_{s=-\infty}^{\infty} \left( \sigma(s) \sigma(s+l-m)+\sigma(s) \sigma(s+l+m) \right)+F \sigma(l) \sigma(m) \right) &\textrm{if}\quad 0<d<\frac{1}{4}
 \end{cases}
\end{eqnarray*}
where $C$, $D$ and $F$ are constants and independent of $l$ and $m$.
Consequently, we can compute:
\begin{displaymath}
\sum_{l=0}^{k}\sum_{m=0}^{k}\frac{\partial g_{i,j} }{\partial x_l \partial x_m}(\sigma(0),\ldots,\sigma(k)) \E((\widehat{\sigma(l)}-\sigma(l))(\widehat{\sigma(m)}-\sigma(m))).
\end{displaymath}
First we study the case $d\geq 1/4$ and we prove that:
\begin{equation}
\sum_{l=0}^k\frac{\partial y_i}{\partial x_l}(\sigma(0),\ldots,\sigma(k))=\left(1-\sum_{r=1}^k a_{r,k} \right) \left( \Sigma_k^{-1}
\left( \begin{array}{c}
1\\
\vdots\\
1
\end{array}\right) \right)_i \label{deriv_d>1/4}
\end{equation}
since if we define $\sigma_0^k:=(\sigma(0),\ldots,\sigma(k))$, we may write the partial derivative as:
\begin{eqnarray}
\frac{\partial y_i}{\partial x_l}(\sigma_0^k)&=&\left( \Sigma_k^{-1}
\left[\frac{\partial 
}{\partial x_l}\left( \begin{array}{c}
x_1 \\
\vdots \\
x_k
\end{array}\right)\right](\sigma_0^k)\right)_i  \nonumber \\ &-&\left(\Sigma_k^{-1}\left[ \frac{\partial 
}{\partial x_l}\left( \begin{array}{cccc}
x_0&x_1& \ldots & x_{k} \\
x_1&x_0& \ddots & x_{k-1} \\
\vdots& \ddots& \ddots& \vdots \\
x_{k}&x_{k-1}&\ldots &x_0
\end{array}\right)\right] (\sigma_0^k)\Sigma_k^{-1} \left( \begin{array}{c}
\sigma(1) \\
\vdots \\
\sigma(k)
\end{array}\right)\right)_i .\label{deriveepartielle}
\end{eqnarray}

Using \eqref{deriv_d>1/4}, the result follows because:
\begin{eqnarray*}
&&\sum_{l=0}^{k}\sum_{m=0}^{k}\frac{\partial g_{i,j} }{\partial x_l \partial x_m}(\sigma(0),\ldots,\sigma(k)) \E((\widehat{\sigma(l)}-\sigma(l))(\widehat{\sigma(m)}-\sigma(m)))\\
&=& \begin{cases}\left(1-\sum_{r=1}^k a_{r,k} \right)^2C n^{4d-2}
\left(\Sigma_k^{-1}
\textbf{1}_{k,k} 
 \Sigma_k^{-1}\right)_{(i,j)}  &\textrm{if} \quad\frac{1}{4}<d<\frac{1}{2} \\
\left(1-\sum_{r=1}^k a_{r,k} \right)^2D n^{-1} \ln(n)
\left(\Sigma_k^{-1}
\textbf{1}_{k,k} 
 \Sigma_k^{-1}\right)_{(i,j)}  &\textrm{if}\quad d=\frac{1}{4} .\end{cases}
\end{eqnarray*}

When $d<1/4$, we first notice by using \eqref{deriveepartielle} that:
\begin{displaymath}
\sum_{l=0}^k\frac{\partial y_i}{\partial x_l}\sigma(l)=\left( -\Sigma_k^{-1} \Sigma_k \Sigma_k^{-1}
\left( \begin{array}{c}
\sigma(1) \\
\vdots\\
\sigma(k)
\end{array}
\right) 
+\Sigma_k^{-1}\left( \begin{array}{c}
\sigma(1) \\
\vdots\\
\sigma(k)
\end{array}
\right) 
\right)_i=0
\end{displaymath}
Then it follows that:
\begin{eqnarray*}
&&\sum_{l=0}^{k}\sum_{m=0}^{k}\frac{\partial g_{i,j} }{\partial x_l \partial x_m}(\sigma(0),\ldots,\sigma(k)) \E((\widehat{\sigma(l)}-\sigma(l))(\widehat{\sigma(m)}-\sigma(m)))\\
&=&\sum_{l=0}^{k}\sum_{m=0}^{k} \frac{\partial y_i}{\partial x_l} \frac{\partial y_j}{\partial x_m} \sum_{s=-\infty}^{\infty} \left( \sigma(s) \sigma(s+l-m) +\sigma(s) \sigma(s+l+m) \right)\\
&=& \frac{1}{2}\sum_{l=0}^{k}\sum_{m=0}^{k} \frac{\partial y_i}{\partial x_l} \frac{\partial y_j}{\partial x_m}\sum_{s=-\infty}^{\infty}\left( \sigma(s) \sigma(s+l-m)+\sigma(s)\sigma(s+m-l) +\sigma(s) \sigma(s+l+m)+\sigma(s)\sigma(s-l-m) \right)\\
&=&\frac{1}{2}\sum_{l=0}^{k}\sum_{m=0}^{k} \frac{\partial y_i}{\partial x_l} \frac{\partial y_j}{\partial x_m}\sum_{s=-\infty}^{\infty}\sigma(s)\int_{-\pi}^{\pi}f(\lambda)\mathrm{e}^{is \lambda}(\mathrm{e}^{i(l-m) \lambda}+\mathrm{e}^{i(m-l) \lambda}+\mathrm{e}^{i(m+l) \lambda}+\mathrm{e}^{i(-m-l) \lambda})\dd \lambda \\
&=& 2\sum_{l=0}^{k}\sum_{m=0}^{k} \frac{\partial y_i}{\partial x_l} \frac{\partial y_j}{\partial x_m}\sum_{s=-\infty}^{\infty}\sigma(s)\int_{-\pi}^{\pi}f(\lambda)\mathrm{e}^{is \lambda}\cos(l \lambda) \cos(m \lambda) \dd \lambda \\
&=& 2 \int_{-\pi}^{\pi}\sum_{s=-\infty}^{\infty}\mathrm{e}^{is \lambda}\sigma(s)f(\lambda) \sum_{l=0}^{k}\sum_{m=0}^{k}\frac{\partial y_i}{\partial x_l} \frac{\partial y_j}{\partial x_m} \cos(l \lambda) \cos(m \lambda) \dd \lambda \\
&=& 2 \int_{-\pi}^{\pi}f(\lambda)^2\sum_{l=0}^{k}\sum_{m=0}^{k}\frac{\partial y_i}{\partial x_l} \frac{\partial y_j}{\partial x_m} \cos(l \lambda) \cos(m \lambda) \dd \lambda \\
&=&2 \left(  \Sigma_k^{-1} H \Sigma_k^{-1}\right) 
\end{eqnarray*}
with $H$ defined in \eqref{H}.
\end{proof}

\subsection{Proof of Lemma \ref{lmmoments}}
\label{lemme2}
Show that:
\begin{equation}
\mathbb{E}\left[ \Big| \frac{1}{n}\sum_{t=1}^{n-k}X_tX_{t+k}-\sigma(k) \Big|^3 \right] =
\begin{cases}
\OO(n^{-3/2})& \textrm{ if $d \leq 1/4$} \\
\OO(n^{6d-3})&\textrm{ if $d > 1/4$}
\end{cases}
\end{equation}
\begin{proof}
\begin{displaymath}
\mathbb{E}\left[ \Big| \frac{1}{n}\sum_{t=1}^{n-k}X_tX_{t+k}-\sigma(k) \Big|^3 \right] \leq \sqrt{\mathbb{E}\left[ \Big| \frac{1}{n}\sum_{t=1}^{n-k}X_tX_{t+k}-\sigma(k) \Big|^2 \right]\mathbb{E}\left[ \Big| \frac{1}{n}\sum_{t=1}^{n-k}X_tX_{t+k}-\sigma(k) \Big|^4 \right]} 
\end{displaymath}
We will separately consider the two terms. First we have:
\begin{eqnarray*}
\mathbb{E}\left[ \Big| \frac{1}{n}\sum_{t=1}^{n-k}X_tX_{t+k}-\sigma(k) \Big|^2 \right]&=&\sigma(k)^2-2\sigma(k) \frac{1}{n}\E\left( \sum_{t=1}^{n-k}X_tX_{t+k}\right)+ \frac{1}{n^2}\E\left( \sum_{t=1}^{n-k}X_tX_{t+k} \sum_{s=1}^{n-k}X_sX_{s+k}\right).
\end{eqnarray*}
Since the process is Gaussian, we have (see \cite{momentsgaussiens}):
\begin{displaymath}
\E \left(X_t X_{t+k}X_s X_{s+k} \right) =\E\left(X_t X_{t+k}\right) \E\left(X_s X_{s+k} \right)+\E\left(X_tX_s \right)\E\left( X_{t+k}X_{s+k}\right) +\E\left( X_t X_{s+k}\right) \E\left(  X_{t+k}X_s\right);
\end{displaymath}
and thus:
\begin{eqnarray*}
\mathbb{E}\left[ \Big| \frac{1}{n}\sum_{t=1}^{n-k}X_tX_{t+k}-\sigma(k) \Big|^2 \right]&=&\left( \frac{(n-k)^2}{n^2}-2\frac{n-k}{n}+1\right)\sigma(k)^2\\&&+\frac{1}{n^2} \sum_{t=1}^{n-k} \sum_{s=1}^{n-k} \sigma(t-s)^2+ \sigma(t+k-s)\sigma(s+k-t) \\
&=&\frac{k^2}{n^2}\sigma(k)^2+\frac{1}{n^2}\sum_{t=1}^{n-k} \sum_{s=1}^{n-k} \left( \sigma(t-s)^2+ \sigma(t+k-s)\sigma(s+k-t)\right) 
\end{eqnarray*}
We note that
\begin{eqnarray*}
\sum_{t=1}^{n-k} \sum_{s=1}^{n-k} \sigma(t-s)^2&=&(n-k)\sigma(0)^2+2\sum_{t=1}^{n-k}(n-k-t)\sigma(t)^2\\
&=& \OO(n) + (n-k)\sum_{t=1}^{n-k}\sigma(t)^2-2\sum_{t=1}^{n-k}t\sigma(t)^2\\
&=& \OO(n)+\OO(n^{4d})
\end{eqnarray*}
In a similar way for:
\begin{displaymath}
\sum_{t=1}^{n-k} \sum_{s=1}^{n-k}\sigma(t+k-s)\sigma(s+k-t)
\end{displaymath}
we obtain that:
\begin{equation}
\sqrt{\mathbb{E}\left[ \Big| \frac{1}{n}\sum_{t=1}^{n-k}X_tX_{t+k}-\sigma(k) \Big|^2 \right]}=
\begin{cases}
\OO(n^{-1/2})& \textrm{ if $d \leq 1/4$} \\
\OO(n^{2d-1}) &\textrm{ if $d > 1/4$}
\end{cases}
\end{equation}
For the second term, we have: 
\begin{eqnarray*}
\mathbb{E}\left[ \Big| \frac{1}{n}\sum_{t=1}^{n-k}X_tX_{t+k}-\sigma(k) \Big|^4 \right]&=&\sigma(k)^4-4\sigma(k)\mathbb{E}\left[ \left( \frac{1}{n}\sum_{t=1}^{n-k}X_tX_{t+k} \right)^3 \right]
+\frac{6\sigma(k)^2}{n^2}\mathbb{E}\left[ \left( \frac{1}{n}\sum_{t=1}^{n-k}X_tX_{t+k} \right)^2 \right]\\
&&-\frac{4\sigma(k)^3}{n}\E\left( \sum_{t=1}^{n-k}X_tX_{t+k}\right) +\mathbb{E}\left[ \left( \frac{1}{n}\sum_{t=1}^{n-k}X_tX_{t+k} \right)^4\right]
\end{eqnarray*}
Since the process is Gaussian, we can apply the result in \cite{momentsgaussiens} and develop the moments as functions which depend only on the covariances of the process. Then we count the order of $\sigma(k)$ in each term of the sum.
The coefficient of $\sigma(k)^4$ is:
\begin{displaymath}
1-\frac{4(n-k)^3}{n^3}+\frac{6(n-k)^2}{n^2}-\frac{4(n-k)}{n}+\frac{(n-k)^4}{n^4}=\frac{k^4}{n^4};
\end{displaymath}
the coefficient of $\sigma(k)^2$ is:
\begin{eqnarray*}
&&\left( \sum_{t=1}^{n-k} \sum_{s=1}^{n-k} \sigma(t-s)^2+ \sigma(t+k-s)\sigma(s+k-t)\right) \left(\frac{-12(n-k)}{n^3}+\frac{6}{n^2}+\frac{6(n-k)^2}{n^4} \right)\\&=&\frac{6k^2}{n^4} \left( \sum_{t=1}^{n-k} \sum_{s=1}^{n-k} \sigma(t-s)^2+ \sigma(t+k-s)\sigma(s+k-t)\right)\\
&=& \begin{cases}
\OO(n^{-3}) & \textrm{ if $d \leq 1/4$} \\
\OO(n^{-4+4d}) &\textrm{ if $d > 1/4$}
    \end{cases}
\end{eqnarray*}
and the coefficient of $\sigma(k)$ is:
\begin{eqnarray*}
&&\left( \frac{1}{n^3}\sum_{t=1}^{n-k}\sum_{s=1}^{n-k}\sum_{r=1}^{n-k}6\sigma(t-s)\sigma(r-s)\sigma(r-t+k)
+\sigma(t+k-r)\sigma(s+k-r)\sigma(r+k-s)\right)\\&\times& \left(\frac{-4}{n^3}+\frac{4(n-k)}{n^4} \right). 
\end{eqnarray*}
We study this asymptotic behaviour as follows:
\begin{eqnarray*}
&&\frac{6}{n^3}\sum_{t=1}^{n-k}\sum_{s=1}^{n-k}\sum_{r=1}^{n-k}\sigma(t-s)\sigma(r-s)\sigma(r-t+k)\\& \leq &\frac{6}{n^3}\sum_{t=1}^{n-k}\sum_{s=1}^{n-k}\sum_{r=1}^{n-k}|\sigma(t-s)\sigma(r-s)\sigma(r-t+k) | \nonumber \\
&\sim&\frac{6}{n^3} \int_{1}^{n-k} \int_{1}^{n-k} \int_{0}^{n-k} |t-s|^{2d-1} |r-s|^{2d-1} |r-t+k|^{2d-1}\dd t \dd s \dd r \nonumber \\
& \leq & \frac{6}{n^3} \int_{1}^{n}\int_{1}^{n} \int_{1}^{n} \int_{1}^{n} |t-s|^{2d-1} |r-s|^{2d-1} |r-t|^{2d-1}\dd t \dd s \dd r \nonumber \\
&\sim& 6n^{6d-3} \int_{0}^{1} \int_{0}^{1} \int_{0}^{1} |t-s|^{2d-1} |r-s|^{2d-1} |r-t|^{2d-1}\dd t \dd s \dd r  \nonumber \\
&=& \mathrm{O}(n^{6d-3}) 
\end{eqnarray*}
The factor of $\sigma(k)$ is bounded by $\OO(n^{6d-4})$.
The constant terms are either like:
\begin{displaymath}
\frac{1}{n^4} \sum_{t=1}^{n-k}\sum_{s=1}^{n-k}\sum_{r=1}^{n-k}\sum_{v=1}^{n-k} \sigma(t-s) \sigma(t-r) \sigma(s-v) \sigma(r-v)
\end{displaymath}
According to a comparison with an integral, they are bounded by $\OO(n^{8d-4})$, or they are like:
\begin{displaymath}
\frac{1}{n^4} \sum_{t=1}^{n-k}\sum_{s=1}^{n-k}\sum_{r=1}^{n-k}\sum_{v=1}^{n-k} \sigma(t-s)^2\sigma(r-v)^2.
\end{displaymath}
We separate the two sums and using the previous results we obtain that:
\begin{equation}
\frac{1}{n^4} \sum_{t=1}^{n-k}\sum_{s=1}^{n-k}\sum_{r=1}^{n-k}\sum_{v=0}^{n-k} \sigma(t-s)^2\sigma(r-v)^2 =
\begin{cases}
\OO(n^{-2}) & \textrm{ if $d \leq 1/4$} \\
\OO(n^{8d-4}) &\textrm{ if $d > 1/4$}
\end{cases}
\end{equation}
When we sum the different components, we obtain that:
\begin{equation}
\sqrt{\mathbb{E}\left[ \Big| \frac{1}{n}\sum_{t=1}^{n-k}X_tX_{t+k}-\sigma(k) \Big|^4\right] }=
\begin{cases}
\OO(n^{-1})& \textrm{ if $d \leq 1/4$} \\
\OO(n^{4d-2})&\textrm{ if $d > 1/4$}
\end{cases}
\end{equation}
Finally, we have obtained that:
\begin{equation}
\mathbb{E}\left[ \Big| \frac{1}{n}\sum_{t=1}^{n-k}X_tX_{t+k}-\sigma(k) \Big|^3 \right] =
\begin{cases}
\OO(n^{-3/2})& \textrm{ if $d \leq 1/4$} \\
\OO(n^{6d-3})&\textrm{ if $d > 1/4$}
\end{cases}
\end{equation}
\end{proof}

\bibliographystyle{apalike} 
\bibliography{biblioprev}

\end{document}